\documentclass[11pt]{article}
\usepackage{amsmath, amsthm, amsfonts, amssymb, color}
\usepackage{mathrsfs}

\usepackage{latexsym,bm}

\usepackage{graphicx}

\newtheorem{thm}{Theorem}[section]

\newtheorem{lemma}[thm]{Lemma}

\newtheorem{remark}[thm]{Remark}
\newtheorem{example}[thm]{Example}
\newtheorem{defi}[thm]{Definition}

\definecolor{wco}{rgb}{0.5,0.2,0.3}

\numberwithin{equation}{section} \theoremstyle{remark}

\newcommand{\gradient}{{\nabla}}

\def\P{{\mathbb P}}
\def\E{{\mathbb E}}
\def\R{{\mathbb R}}
\def\1{{\mathbbm 1}}
\def\eps{{\varepsilon}}

\def\LL{{\mathcal L}}
\def\SS{{\mathcal S}}

\numberwithin{equation}{section} 

\begin{document}

\allowdisplaybreaks

\title{\bf Boundary Harnack principle and gradient estimates for fractional Laplacian perturbed by non-local operators}
\author{{\bf Zhen-Qing Chen\footnote{Department of Mathematics, University of Washington, Seattle,
WA 98195, USA. E-mail: zqchen@uw.edu}},
\quad {\bf Yan-Xia Ren}\footnote{LMAM School of Mathematical Sciences \& Center for
Statistical Science, Peking University, Beijing, 100871, P.R. China. E-mail: yxren@math.pku.edu.cn} \quad  \hbox{and} \quad {\bf Ting
Yang}\footnote{Corresponding author. School of Mathematics and Statistics, Beijing Institute of Technology, Beijing 100081, P.R. China. Email: yangt@bit.edu.cn}
\footnote{Beijing Key Laboratory on MCAACI, Beijing 100081, P.R. China.}
}
\date{}
\maketitle

\begin{abstract}
Suppose $d\ge 2$ and $0<\beta<\alpha<2$. We consider the non-local
operator $\mathcal{L}^{b}=\Delta^{\alpha/2}+\mathcal{S}^{b}$, where
$$
\mathcal{S}^{b}f(x):=\lim_{\varepsilon\to
0}\mathcal{A}(d,-\beta)\int_{|z|>\varepsilon}\left(f(x+z)-f(x)\right)
\frac{b(x,z)}{|z|^{d+\beta}}\,dy.
$$
Here $b(x,z)$ is a bounded measurable function on
$\mathbb{R}^{d}\times\mathbb{R}^{d}$
that is symmetric in $z$, and $\mathcal{A}(d,-\beta)$ is a normalizing constant so that when $b(x, z)\equiv 1$, $\mathcal{S}^{b}$
becomes
 the fractional Laplacian
$\Delta^{\beta/2}:=-(-\Delta)^{\beta/2}$. In other words,
$$
 \mathcal{L}^{b}f(x):=\lim_{\varepsilon\to
0}\mathcal{A}(d,-\beta)\int_{|z|>\varepsilon}\left(f(x+z)-f(x)\right) j^b(x, z)\,dz,
$$
where
$j^b(x, z):= \mathcal{A}(d,-\alpha) |z|^{-(d+\alpha)}
+ \mathcal{A}(d,-\beta) b(x, z)|z|^{-(d+\beta)}$.
It is recently established in Chen and Wang \cite{CW} that,
when $j^b(x, z)\geq 0$ on $\R^d\times \R^d$, there is a
conservative  Feller process $X^{b}$ having $\LL^b$ as its infinitesimal
generator. In this paper we establish,
under certain conditions on $b$,
a uniform boundary Harnack principle for harmonic functions of $X^b$ (or equivalently, of $\LL^b$) in
any $\kappa$-fat open set.
We further establish uniform gradient estimates for
non-negative harmonic functions of $X^{b}$  in open sets.
\end{abstract}

\medskip


\medskip


\section{Introduction}

Let  $d\ge 2$,  $0<\beta<\alpha<2$,
and  $b(x,z)$ be a bounded measurable function on
$\mathbb{R}^{d}\times\mathbb{R}^{d}$ with $b(x,z)=b(x,-z)$ for
$x,z\in \mathbb{R}^{d}$.
Consider the non-local
operator $\mathcal{L}^{b}=\Delta^{\alpha/2}+\mathcal{S}^{b}$, where
\begin{equation}
\mathcal{S}^{b}f(x):=\lim_{\varepsilon\to
0}\mathcal{A}(d,-\beta)\int_{|z|>\varepsilon}\left(f(x+z)-f(x)\right)\frac{b(x,z)}{|z|^{d+\beta}}\,dz.
\nonumber
\end{equation}
Here $\mathcal{A}(d,-\beta)$ is a normalizing constant so that when $b(x, z)\equiv 1$,
$ \mathcal{S}^{b}$ becomes the fractional Laplacian $ \Delta^{\beta/2}:=-(-\Delta)^{\beta/2}$;
in other words,
$\mathcal{A}(d,-\beta)=\beta 2^{\beta-1}\pi^{-d/2}\Gamma((d+\beta)/2)/ \Gamma(1-\beta/2)$.
Thus $\LL^b$ can be expressed as
\begin{equation}
\mathcal{L}^{b}f(x)=\lim_{\varepsilon\to 0}\int_{|z|>\varepsilon}\left(f(x+z)-f(x)\right) j^b(x, z)\,dz
\nonumber
\end{equation}
where
$$
j^b (x, z)= \frac{\mathcal{A}(d,-\alpha)}{|z|^{d+\alpha}}+
\frac{\mathcal{A}(d,-\beta) b(x, z)}{|z|^{d+\beta}}.
$$
Note that since $b(x, z)$ is symmetric in $z$, for $f\in C^2_b(\R^d)$,
$$
\mathcal{L}^{b}f(x) = \int_{\R^d} \left(f(x+z)-f(x)-\nabla f(x)
\cdot z \mathbf{1}_{\{|z|\leq 1\}}\right) j^b(x, z)\,dz.
$$
Recently, $ \mathcal{L}^{b}$,
the fractional Laplacian perturbed by a lower order non-local operator $\SS^b$,
and its fundamental solution have been studied in Chen and Wang \cite{CW}.
It is established there that if for
every $x\in\mathbb{R}^{d}$,
$j^b(x, z)\geq 0$  (that is, $b(x,z)\ge -\frac{\mathcal{A}(d,-\alpha)}{\mathcal{A}(d,-\beta)}\,|z|^{\beta-\alpha}$)
for a.e. $z\in\mathbb{R}^{d}$, then $\mathcal{L}^{b}$ has a unique jointly continuous
fundamental solution $p^{b}(t,x,y)$,  which
uniquely determines a conservative Feller process $X^{b}$ on
the canonical Skorokhod space
$\mathbb{D}([0,+\infty),\mathbb{R}^{d})$ such that
$$\mathbb{E}_{x}\left[f(X^{b}_{t})\right]=\int_{\mathbb{R}^{d}}f(y)p^{b}(t,x,y)dy,
\qquad x\in \mathbb{R}^{d}.
$$
for every bounded measurable function $f$ on $\R^d$.
The Feller process $X^{b}$ is typically non-symmetric and it has a L\'evy
system
  $(J^{b}(x,y)dy,t)$
  (see \cite[Theorem 1.3]{CW}),
  where
\begin{equation}\label{e:1.5}
J^{b}(x,y):=j^b(x, y-x).
\end{equation}
When $b$ takes constant
value $\eps >0$,
$X^\eps$  has the same distribution as the L\'evy process $Y+\eps^{1/\beta} Z$, where $Y$ and $Z$
are rotationally symmetric $\alpha$- and $\beta$-stable processes
on $\R^d$ that are independent of each other.
Moreover, two-sided heat kernel estimates  for  $\mathcal{L}^{b}$
have been obtained in Chen and Wang \cite{CW},
 while two-sided Dirichlet kernel estimates in $C^{1,1}$
open sets have recently been
obtained in Chen and Yang \cite{CY}.
We emphasize that the lower order perturbations $ \mathcal{S}^{b}$ considered in \cite{CW, CY} and in this paper
are of high intensity, in the sense that $\int_{\R^d} \frac{|b(x, z)|}{|z|^{d+\beta}} dz$ may be identically infinite.
This is the case, for example, when $b(x, z)$ is bounded away from $0$.
Such high intensity non-local perturbation is in stark contrast to the Feynman-Kac type non-local perturbation considered
in  \cite{BS,CKS4, CS2, GR,KL}, where the kernels used in the non-local perturbation are typically  integrable.
Under high intensity non-local perturbation, the potential theory for the perturbed generator $ \mathcal{L}^{b}$
can be different from that of $\Delta^{\alpha/2}$. For example, the heat kernel of  $ \mathcal{L}^{b}$
may not be comparable to that of $\Delta^{\alpha/2}$; see \cite{CW} for details.

In this paper, we investigate boundary Harnack principle and gradient estimates
for non-negative harmonic functions of $\mathcal{L}^{b}$ in open sets.
Boundary Harnack principle (BHP) asserts
that non-negative harmonic functions that vanish
in an exterior part of a neighborhood at the boundary decay at the same rate.
It is an important property in analysis and in probability theory on harmonic functions.
We refer the reader to the introduction of \cite{BKK, KSV2}
for a brief account on the history of BHP
that started with Brownian motion
and then extended to
subordinate Brownian
motions
and to certain pure jump strong Markov processes. Since
$\mathcal{L}^{b}$ is typically
state-dependent
and its dual operator may not be
Markovian,
the BHP results in \cite{BKK, KSV2} are not applicable to harmonic functions
of $ \mathcal{L}^{b}$.
In this paper, we establish uniform boundary Harnack principle
on $\kappa$-fat open sets for
non-negative harmonic functions
of  $\mathcal{L}^{b}$
by estimating Poisson kernels of $\LL^b$ in small balls.

Gradient estimates for harmonic functions of elliptic operators and on manifolds
have been studied extensively in literature, including the celebrated Li-Yau inequality.
See \cite{Cr} and the references therein and for a coupling argument.
See also \cite{SSW} for  gradient estimates of the transition semigroups of some L\'{e}vy processes
using coupling method.
Gradient estimates for harmonic functions for non-local operators are quite recently.
In \cite{BKN}, a gradient estimate for harmonic functions of symmetric stable processes
is obtained. Gradient estimates for harmonic functions of mixed stable processes
were derived in \cite{Y}.
It has recently been extended to a class of isotropic unimodal  L\'evy process in \cite{KR}. For gradient estimate for harmonic functions of the Schr\"odinger operator
$\Delta^{\alpha/2}+q$, see \cite{BKN} for $\alpha \in (1, 2)$
and \cite{Ku} for $\alpha\in (0, 1]$.
The second main result of this paper is to establish gradient estimates
for positive harmonic functions of  $\mathcal{L}^{b}$. As far as we know, this is
the first gradient estimate result for non-L\'evy non-local operators.

We now describe our main results in details.
In this paper, we use ``$:=$" as a way of definition.  For $a,b\in \mathbb{R}$, $a\wedge
b:=\min\{a,b\}$ and $a\vee b:=\max\{a,b\}$.  Let $|x-y|$ denote the Euclidean distance
between $x$ and $y$, and  $B(x,r)$ the open ball centered at $x$
with radius $r>0$.
For any two positive
functions $f$ and $g$, $f\stackrel{c}\lesssim g$ means that there is
a positive constant $c$ such that $f\le c g$ on their common domain
of definition, and $f\stackrel{c}\asymp g$ means that $c^{-1}g\le
f\le c g$. We also write ``$\lesssim$" and ``$\asymp$" if $c$ is
unimportant or understood.
If $D\subset \mathbb{R}^{d}$ is an open set, for
every $x,y\in D$, define
\begin{equation}\label{e:1.3}
\delta_{D}(x):=\mathrm{dist}(x,\partial D) \quad \hbox{and} \quad  r_{D}(x,y):=\delta_{D}(x)+\delta_{D}(y)+|x-y|.
\end{equation}
It is easy
to see that
\begin{equation}\label{e:1.4}
r_{D}(x,y)\asymp \delta_{D}(x)+|x-y|\asymp \delta_{D}(y)+|x-y|.
\end{equation}
Denote by $\tau^{b}_{D}:=\inf\{t>0:X^{b}_{t}\not\in D\}$, the exit time from $D$ by $X^b$.
When there is no danger of confusion, we will drop the superscript $b$ and simply write $\tau_D$ for $\tau^b_D$.

\begin{defi}\label{D:1}  \rm
A function $f$ defined on $\mathbb{R}^{d}$ is said to be harmonic in an open
set $D$ with respect to $X^{b}$ if it has the mean-value property:
for every bounded open set $U\subset D$ with $\overline{U}\subset
D$,
\begin{equation}\label{harmonic}
f(x)=\mathbb{E}_{x}\left[ f(X^{b}_{\tau_{U}})\right]\quad\hbox{for }
x\in U.
\end{equation}
It is said to be  regular harmonic in $D$ if
\eqref{harmonic} holds for $U=D$.
\end{defi}

Denote by  $\partial$  a cemetery point that is added to $D$ as an isolated point.
We use the convention that $X^b_\infty:=\partial$
and any function $f$ is extended to the cemetery point $\partial $ by setting
$f(\partial )=0$.
So $\mathbb{E}_{x}\left[f(X^{b}_{\tau_{D}})\right]$ should be
understood as
$\mathbb{E}_{x}\left[f(X^{b}_{\tau_{D}}):\tau_{D}<+\infty\right]$.
In Definition \ref{D:1}, we always assume implicitly that the
expectation in \eqref{harmonic} is absolutely convergent.

\medskip

\noindent\textbf{Assumption 1}
Suppose $M_1, M_2 \geq 1$.
 $b(x, z)$ is a bounded function
on $\R^d\times \R^d$ satisfying
\begin{equation}\label{condi0}
\|b\|_{\infty}\le {M_1} \quad \hbox{ and } \quad
b(x,z)=b(x,-z) \quad \hbox{for }  x,z\in\mathbb{R}^{d},
\end{equation}
and there exists a positive constant $\eps_0 \in [0, 1]$ such that for
every $x,y\in\mathbb{R}^{d}$,
\begin{equation} \label{condi2}
M^{-1}_{2}J^{\eps_0}(x,y)\le J^{b}(x,y)\le {M_2}
J^{\eps_0}(x,y).
\end{equation}
Here $J^{\eps_{0}}(x,y)=\mathcal{A}(d,-\alpha)|x-y|^{-d-\alpha}+\eps_{0}\mathcal{A}(d,-\beta)|x-y|^{-d-\beta}$.
Since $J^{\eps_{0}}(x,y)$ depends only on $|x-y|$, we also write $J^{\eps_{0}}(|x-y|)$ for $J^{\eps_{0}}(x,y)$.

\begin{defi}\rm Let $\kappa \in (0, 1)$. An open set $D\subset \R^d$ is said to be
$\kappa$-fat
if for every $z\in \partial D$ and $r\in (0, 1]$, there is some point
$x\in D$ so that $B(x, \kappa r) \subset D\cap B(z, r)$.
\end{defi}

The following is the first main result of this paper.

\begin{thm}[Uniform boundary Harnack inequality]\label{BHP}
Suppose  Assumption 1 holds and $D$ is a $\kappa$-fat open set in $\R^d$ with $\kappa \in (0, 1)$.
There exist constants $r_{1}=r_{1}(d,\alpha,\beta,{M_1})\in (0,1]$ and $C_{1}=C_{1}(d,\alpha,\beta,\kappa,{M_1},{M_2})\ge 1$
such that for
every $z_0\in \partial D$ and $r\in (0,r_{1}/2]$, and all non-negative
functions $u,v$ that are regular harmonic in $D\cap B(z_0,2r)$ with
respect to $X^{b}$ and vanish in $D^{c}\cap B(z_0,2r)$, we have
$$
 \frac{u(x)}{v(x)}\le C_{1}\frac{u(y)}{v(y)}
 \qquad \hbox{for }  x, y\in D\cap B(z_0,r).$$
\end{thm}

We call the above property {\it uniform} boundary Harnack principle  because
the constants $r_1$ and $C_1$ in the above theorem
are independent of $\eps_0 \in [0, 1]$
appeared in
condition \eqref{condi2}.
We next study the gradient estimates for non-negative harmonic functions in open sets.
We write $\partial_{x_i}$ or $\partial_i$ for $\frac{\partial}{\partial
x_{i}}$ and $\gradient$ for $( \partial_{x_1} ,\cdots, \partial_{x_d})$.

\begin{thm}\label{theorem1}
Let $D$ be an arbitrary open set in $\mathbb{R}^{d}$.
Under Assumption 1, there is a
constant $C_{2}=C_{2}(d,\alpha,\beta,{M_1},{M_2})>0$ such that for
any non-negative function $f$ in $\mathbb{R}^{d}$ which is harmonic in
$D$ with respect to $X^{b}$, $\nabla f(x)$ exists for every $x\in
D$, and we have
\begin{equation}
|\nabla f(x)|\le C_{2}\frac{f(x)}{1\wedge \delta_{D}(x)}
\quad\hbox{for }  x\in D.
\nonumber
\end{equation}
\end{thm}

For $x=(x_{1},\cdots,x_{d})\in\mathbb{R}^{d}$ and $1\leq i\leq d$, we write
$\tilde{x}^{i}$ for $(x_{1},\cdots,x_{i-1},x_{i+1},\cdots,x_{d})\in\mathbb{R}^{d-1}$.

\medskip

\noindent \textbf{Assumption 2.}
Suppose there is  $i\in\{1,\cdots,d\}$
so that for every $x\in\mathbb{R}^{d}$,
\begin{equation}
b(x,z)=\varphi(\tilde{x}^{i})\psi(|z|)\quad\mbox{a.e.
}z\in\mathbb{R}^{d},
\nonumber
\end{equation}
where $\varphi:\mathbb{R}^{d-1}\to \mathbb{R}$ is a non-negative
measurable function, and $\psi:\mathbb{R}_{+}\to \mathbb{R}$ is
a measurable function such that
\begin{equation}\label{e:1.10}
\frac{\psi ( r)}{r^{d+\beta}} \
\hbox{ is non-increasing  in } r>0.
\end{equation}

\begin{thm}\label{theorem2}
Suppose  Assumption 1 and Assumption 2 hold.  Let
$D=\left\{ x\in\mathbb{R}^{d}:\ x_{i}>\Gamma(\tilde{x}^{i}) \right\}$
be a special Lipschitz domain, where
$\Gamma:\mathbb{R}^{d-1}\to\mathbb{R}$ is a Lipschitz
function with Lipschitz constant $\lambda_{0}$ (that is,
$|\Gamma(\tilde{x}^{i})-\Gamma(\tilde{y}^{i})|\le\lambda_{0}|\tilde{x}^{i}-\tilde{y}^{i}|$
for every $\tilde{x}^{i},\tilde{y}^{i}\in\mathbb{R}^{d-1}$).
Then there are positive constants
$R_{1}=R_{1}(d,\alpha,\beta,\lambda_{0},{M_1}, {M_2})$ and
$C_{3}=C_{3}(d,\alpha,\beta,\lambda_{0},{M_1},{M_2})\ge 1$
such that for every $r\in(0,R_{1}]$,
there is a constant $\eta_{1}=\eta_{1}(d,\alpha,\beta,\lambda_{0},{M_1},{M_2},r)\in (0,r/2)$ so that for every $z_0\in\partial D$  and every
non-negative function $f$ that is harmonic in $D\cap B(z_0,r)$ with
respect to $X^{b}$ and vanishes in $D^{c}\cap B(z_0,r)$,
\begin{equation}\label{them2.main}
C_{3}^{-1}\frac{f(x)}{\delta_{D}(x)} \le |\nabla f(x)|\le C_{3}
\frac{f(x)}{\delta_{D}(x)} \quad \hbox{for }
x\in D\cap B(z_{0},\eta_{1}).
\end{equation}
\end{thm}

Obviously Assumption 2 is implied by

\noindent \textbf{Assumption 3.} There exists a measurable function
$\psi:\mathbb{R}_{+}\to \mathbb{R}$ satisfying \eqref{e:1.10}
such that for every $x\in\mathbb{R}^{d}$,
\begin{equation}
b(x,z)=\psi(|z|)\quad\mbox{a.e. }z\in\mathbb{R}^{d}.
\nonumber
\end{equation}

\begin{defi}\rm An open set $D\subset \mathbb{R}^{d}$ is said to be
 Lipschitz if  for every $z_0\in \partial D$,
there is a Lipschitz function
$\Gamma_{z_0}:\mathbb{R}^{d-1}\rightarrow \mathbb{R}$, an orthonormal
coordinate system $CS_{z_0}$ and
a constant $R_{z_0}>0$
such that if
$y=(y_{1},\cdots,y_{d-1},y_{d})$ in $CS_{z_0}$ coordinates, then
$$D\cap B(z_0,R_{z_0})=\{y:y_{d}>\Gamma_{z_0}(y_{1},\cdots,y_{d-1})\}\cap B(z_0,R_{z_0}).$$
If there exist positive constants $R_{0}$ and $\lambda_0$ so that
$R_{z_0}$ can be taken to be $R_0$ for all $z_0\in \partial D$ and the
Lipschitz constants of $\Gamma_{z_0}$ are not greater than
$\lambda_{0}$, we call $D$ a Lipschitz open set with characteristics
$(\lambda_{0},R_{0})$.
\end{defi}

Clearly, if $D$ is a Lipschitz open set with characteristics $(\lambda_{0},R_{0})$,
then it is $\kappa$-fat for some $\kappa =\kappa(\lambda_{0},R_{0})\in (0,1)$.
The following theorem follows directly from Theorem \ref{theorem2}.

\begin{thm}\label{T:1.7}
Let $D$ be a Lipschitz open set in $\mathbb{R}^{d}$ with
characteristics $(\lambda_0, R_0)$. Under Assumptions 1 and 3,
there are positive  constants
$R_{2}=R_{2}(d,\alpha,\beta,\lambda_{0},R_{0},{M_1}, {M_2})$
and \hfill \break
$C_{4}=C_{4}(d,\alpha,\beta,\lambda_{0},R_{0},{M_1},{M_2}) \geq 1$
such that for  every $r\in(0,R_{2}]$,
there is a constant $\eta_{2}=\eta_{2}(d,\alpha,\beta,\lambda_{0},R_{0},{M_1},{M_2},r)\in (0,r/2)$ so that for every  $z_0\in\partial D$  and
every non-negative function $f$ that is harmonic in $D\cap B(z_0,r)$
with respect to $X^{b}$ and vanishes in $D^{c}\cap B(z_0,r)$,
\begin{equation}
C_{4}^{-1}\frac{f(x)}{\delta_{D}(x)} \le |\nabla f(x)|\le C_{4}
\frac{f(x)}{\delta_{D}(x)} \qquad\hbox{for }
x\in D\cap B(z_{0},\eta_{2}).
\nonumber
\end{equation}
\end{thm}

Results in Theorem \ref{theorem1}, Theorem \ref{theorem2} and Theorem \ref{T:1.7}
can be called {\it uniform} gradient estimates because the constants $C_k$, $2\leq k\leq 4$, and $\eta_i$, $1\leq i\leq 2$, are independent of $\eps_0$ of \eqref{condi2}.

The rest of the paper is organized as follows.
Preliminary results on Green functions and Poisson kernels are presented in
Section \ref{S:2}. The proof of the uniform boundary Harnack principle
is given in Section \ref{S:3}.
Section \ref{S:4} is devoted to the proof of Theorem \ref{theorem1},
while the proof of Theorem \ref{theorem2} is given in Section \ref{S:5}.
In this paper,
we  use capital letters
$C_{1}, C_{2}, \cdots$
to denote constants in the
statements of results. The lower case constants
$c_{1},c_{2},\cdots,$ will denote the generic constants used in
proofs, whose exact values are not important, and can change from
one appearance to another. We use $e_k$ to denote the unit vector along
the positive direct of $x_k$-axis.

\section{Preliminaries}\label{S:2}

Recall the  L\'{e}vy system  $(J^{b}(x,y)dy,t)$ from \eqref{e:1.5}, which describes the jumps of $X^{b}$: for
any non-negative measurable function $f$ on
$\mathbb{R}_{+}\times\mathbb{R}^{d}\times\mathbb{R}^{d}$ with
$f(s,y,y)=0$ for all $y\in \mathbb{R}^{d},\ x\in \mathbb{R}^{d}$ and
stopping time $T$ (with respect to the filtration of $X^{b}$),
\begin{equation}
\mathbb{E}_{x}\left[\sum_{s\le
T}f(s,X^{b}_{s-},X^{b}_{s})\right]=\mathbb{E}_{x}\left[\int_{0}^{T}\int_{\mathbb{R}^{d}}f(s,X^{b}_{s},y)J^{b}(X^{b}_{s},y)dy\,ds\right].\label{levysystem}
\end{equation}
Suppose $D$ is a
Greenian open set of $X^{b}$. Let $G^{b}_{D}(x,y)$ denote the Green
function of $D$, that is, $$
\int_{D}f(y)G^{b}_{D}(x,y)dy=\mathbb{E}_{x}\left[\int_{0}^{\tau_{D}}f(X^{b}_{s})ds\right]
$$
for every bounded measurable function $f$ on $D$ and $x\in D$.
It follows from \eqref{levysystem} that for every bounded open set
$D$ in $\mathbb{R}^{d}$, every $f\ge 0$, and $x\in D$,
\begin{equation}
\mathbb{E}_{x} \left[ f(X^{b}_{\tau_{D}}):X^{b}_{\tau_{D}-}\not=X^{b}_{\tau_{D}}\right]
=\int_{\bar{D}^{c}}f(z)\left(\int_{D}G^{b}_{D}(x,y)J^{b}(y,z)dy\right)dz.\label{(*)}
\end{equation}
Define
\begin{equation}\label{(3)}
K^{b}_{D}(x,z):=\int_{D}G^{b}_{D}(x,y)J^{b}(y,z)dy
\quad\hbox{for } (x,z)\in D\times \bar{D}^{c}.
\end{equation}
We call $K^b_D(x, z)$
the Poisson kernel of $X^b$ on $D$.
Then \eqref{(*)} can be written as
\begin{equation}
\mathbb{E}_{x}\left[f(X^{b}_{\tau_{D}}):X^{b}_{\tau_{D}-}\not=X^{b}_{\tau_{D}}\right]
=\int_{\bar{D}^{c}} f(z) K^{b}_{D}(x,z) dz.
\nonumber
\end{equation}
For any $\lambda>0$, define
\begin{equation}
b_{\lambda}(x,z):=\lambda^{\beta-\alpha}b(\lambda^{-1}x,\lambda^{-1}z)  \quad
\hbox{for }  x,z\in\mathbb{R}^{d}.
\nonumber
\end{equation}
It is not hard to show that
\begin{equation} \label{e:2.6a}
J^{b_{\lambda}} (x, y) = \lambda^{-(d+\alpha)} J^b (\lambda^{-1}x, \lambda^{-1}y)\quad \hbox{for }  x,y\in\mathbb{R}^{d}.
\end{equation}
and
\begin{equation}
\{\lambda X^{b}_{\lambda^{-\alpha}t}; t\geq 0\} \hbox{
has the same distribution as } \{ X^{b_{\lambda}}_{t}; t\geq 0\}.
\nonumber
\end{equation}
 So for any $\lambda>0$, we have the following scaling properties:
\begin{equation}\label{scalingforgb}
G^{b}_{D}(x,y)=\lambda^{d-\alpha}G^{b_{\lambda}}_{\lambda D}(\lambda
x,\lambda y)\quad\hbox{for }  x,y\in D,
\end{equation}
\begin{equation}
K^{b}_{D}(x,z)=\lambda^{d}K^{b_{\lambda}}_{\lambda D}(\lambda
x,\lambda z)\quad\hbox{for }  x\in D,\
z\in\bar{D}^{c}.\label{scallingpoisson}
\end{equation}
 If $u$ is
harmonic in $D$ with respect to $X^{b}$, then for any $\lambda>0$,
$v(x):=u(x/\lambda)$ is harmonic in $\lambda D$ with respect to
$X^{b_{\lambda}}$.

When $b(x, z)\equiv 0$, $X^0$ is simply
an isotropic  symmetric $\alpha$-stable
process on $\mathbb{R}^{d}$, which we will denote as $X$.
We will also write $J$ for $J^0$.
 It is known that if $d>\alpha$,  the
process $X$ is transient and its Green function  is given by
\begin{equation}\label{globalgreen}
G(x,y)=\frac{\Gamma(d/2)} {2^{\alpha}\pi^{d/2}
\Gamma(\alpha/2)^{2}}\, |x-y|^{\alpha-d}
\quad\hbox{for }  x,y\in\mathbb{R}^{d}.
\end{equation}
It is shown in Blumenthal \textit{et al.} \cite{Blumenthal}
that the Green function  of $X$ in a ball $B(0, r)$ is given by
\begin{equation}\label{ballgreen}
G_{B(0,r)}(x,y)=\frac{\Gamma(d/2)} {2^{\alpha}\pi^{d/2}
\Gamma(\alpha/2)^{2}}
\int_{0}^{z}(u+1)^{-d/2}u^{\alpha/2-1}du\,|x-y|^{\alpha-d}\quad\hbox{for }
x,y\in B(0,r),
\end{equation}
where $z=(r^{2}-|x|^{2})(r^{2}-|y|^{2})|x-y|^{-2}$ and $r>0$.
The above formula yields the following two-sided estimates (see, for example,
\cite{Chen}):
Suppose $B$ is an arbitrary ball in $\mathbb{R}^{d}$ with radius
$r>0$. Then there is a universal  constant $c_{1}=c_{1}(d,\alpha)>1$
so that for every $x,y\in B$,
\begin{equation}\label{e:2.10}
  G_{B}(x,y) \stackrel{c_{1}}\asymp|x-y|^{\alpha -d}
  \left(1\wedge \frac{\delta_{B}(x) }{|x-y|}\right)^{\alpha/2}
  \left(1\wedge \frac{\delta_{B}(y) }{|x-y|}\right)^{\alpha/2}.
\end{equation}
Since for $a, b>0$, $a\wedge b \asymp \frac{ab}{a+b}$ and
$1\wedge \frac{a}{b}\asymp \frac{a}{a+b}$, in view of \eqref{e:1.4}
we can rewrite \eqref{e:2.10} as
\begin{equation}\label{e:2.11}
  G_{B}(x,y) \asymp|x-y|^{\alpha -d}
    \frac{\delta_{B}(x)^{\alpha/2} \delta_{B}(y)^{\alpha/2}}
  {r_B (x, y)^\alpha}.
\end{equation}
It follows immediately from \eqref{e:2.10} that
there is a positive constant
$c_{2}=c_{2}(d,\alpha)>1$ so that for $B=B(x_0, r)$,
$$
c_{2}^{-1}r^{\alpha}\le \mathbb{E}_{x}\tau_{B}\le c_{2}r^{\alpha}
\qquad  \hbox{for } x\in B(x_0, r/2).
$$
 Riesz (see \cite{Blumenthal}) derived the
following explicit formula for  the Poisson kernel $K_{B(0,r)}(x,z)$ of $X$ on
$B(0, r)$.
\begin{equation}\label{ballpoisson}
K_{B(0,r)}(x,z)=  \frac{\Gamma(d/2)\sin (\pi\alpha / 2) }{\pi^{d/2+1}} \, \frac{(r^{2}-|x|^{2})^{\alpha/2}}{
(|z|^{2}-r^{2})^{\alpha/2} \, |x-z|^{d}} \quad\hbox{for }  |x|<r \hbox{ and }
 |z|>r ,
\end{equation}
We also know from \cite [Lemma 6]{Bogdan1} that, for rotationally symmetric $\alpha$-stable process $X$,
$$
 \P_x (X_{\tau_D} \in \partial D) =0 \quad \hbox{and so} \quad \P_x (X_{\tau_D}\not= X_{\tau_D -})= 1
$$
for every $x\in D$ if $D$ is a domain that satisfies uniform exterior cone condition.
(Although the result of  \cite [Lemma 6]{Bogdan1} is stated only for bounded domain satisfying uniform exterior
cone condition, its proof works for any domain satisfying uniform exterior
cone condition.)

\section{Boundary Harnack principle}\label{S:3}

Recall that we write $X$ and $J$ for $X^0$ and $J^0$, respectively.
First we record the following gradient estimate on the Green function $G_D$
of symmetric $\alpha$-stable process $X$ from \cite{BKN}.

\begin{lemma}\label{lemma1}
{\rm (\cite[ Corollary 3.3]{BKN})}
 Let $D$ be a Greenian
domain in $\mathbb{R}^{d}$ of $X$. Then
\begin{equation}
|\nabla G_{D}(x,y)|\le
d\frac{G_{D}(x,y)}{|x-y|\wedge\delta_{D}(x)} \quad \hbox{for }  x,y\in
D,\ x\neq y.\label{partialgreen}
\end{equation}
\end{lemma}

For $x\not= y$ in $D$, define
\begin{equation} \label{e:3.2}
h_{D}(x,y):=\begin{cases}
         |x-y|^{\alpha-\beta-d}\left(1\wedge\frac{\delta_{D}(y)}{|x-y|}\right)^{\alpha/2}
         \quad &\hbox{if } \alpha >2 \beta , \smallskip \\
         |x-y|^{\beta-d}\left(1\wedge\frac{\delta_{D}(y)}{|x-y|}\right)^{\beta}
         \left(1\vee\log\frac{|x-y|}{\delta_{D}(x)}\right)\quad &\hbox{if } \alpha =2\beta ,
         \smallskip \\
         |x-y|^{\alpha-\beta-d}\left(1\wedge\frac{\delta_{D}(y)}{|x-y|}\right)^{\alpha/2}
         \left(1\vee\frac{|x-y|}{\delta_{D}(x)}\right)^{\beta-\alpha/2}
         &\hbox{if }   \alpha<2 \beta ,
        \end{cases}
\end{equation}

The following two results are established in \cite{CY}.

\begin{lemma} {\rm (\cite[ Theorems 4.11 and Lemma 4.13]{CY})} \label{lemma2}
Suppose $b$ is a bounded function satisfying \eqref{condi0},
  and that for every $x\in\mathbb{R}^{d}$,
$J^{b}(x,y)\ge 0$ a.e. $y\in \mathbb{R}^{d}$. There exist
positive constants $r_{1}=r_{1}(d,\alpha,\beta,{M_1})\in (0, 1]$
and $C_{5}=C_{5} (d, \alpha, \beta, {M_1})$ such that for any
$x_{0}\in\mathbb{R}^{d}$ and any ball $B=B(x_{0},r)$ with radius
$r\in (0,r_{1}]$, we have for  $x,y\in B$,
\begin{equation} \label{e:3.3}
\frac{1}{2}G_{B}(x,y)\le G^{b}_{B}(x,y)\le
\frac{3}{2}G_{B}(x,y) \quad \hbox{ and } \quad
|\mathcal{S}^{b}_{x}G^{b}_{B}(x,y)|\le C_{5}  h_{B}(x,y)
\end{equation}
Moreover, $\mathbb{P}_{x}\left(X^{b}_{\tau_{B}}\in\partial
B\right)=0$ for every $x\in B$. In this case, for every
non-negative measurable function $f$,
$$
\mathbb{E}_{x}f(X^{b}_{\tau_{B}})=\int_{\bar{B}^{c}}f(z)K^{b}_{B}(x,z)dz\quad\hbox{for }  x\in B.
$$
\end{lemma}

\begin{lemma} {\rm (\cite[Lemma 4.2]{CY})} \label{lemma10}
Let $D$ be a bounded open set in $\mathbb{R}^{d}$. There exists a
constant $C_{6}=C_{6}(d,\alpha,\beta,\mathrm{diam}(D),{M_1})>0$ such
that for any bounded function $b$ satisfying \eqref{condi0},
 and that for every $x\in\mathbb{R}^{d}$,
$J^{b}(x,y)\ge 0$ for a.e. $y\in\mathbb{R}^{d}$, we have
$$G^{b}_{D}(x,y)\le C_{6}|x-y|^{\alpha -d}\quad\hbox{for }  x,y\in D.$$
\end{lemma}

Note that the constant $C_7$ below is independent of $\eps_0\in [0, 1]$
appeared in \eqref{condi2}.

\begin{thm}[Uniform Harnack inequality]\label{HI1}
Let $r_1\in (0, 1]$ be the constant in Lemma \ref{lemma2}.
Under Assumption 1,
there exists a constant $C_{7}=C_{7}(d,\alpha,\beta,{M_1},{M_2})\ge
1$ such that for  every $x_{0}\in
\mathbb{R}^{d}$, $r\in (0,r_{1}]$, and every non-negative function $u$
which is regular harmonic in $B(x_{0},r)$, we have
$$\sup_{y\in B(x_{0},r/2)}u(y)\le C_{7}\inf_{y\in B(x_{0},r/2)}u(y).$$
\end{thm}

\proof   Let
$u^{*}(x):=\mathbb{E}_{x}\left[u(X^{\eps_0}_{\tau_{B(x_{0},r)}})\right]$.
Then $u^{*}$ is regular harmonic in
$B(x_{0},r)$ with respect to the mixed stable processes $X^{\eps_0}$.
In view of Lemma \ref{lemma2} and Assumption 1,
for every $x_{0}\in \mathbb{R}^{d}$ and $r\in (0,r_{1}]$,
the Poisson kernel $K^{b}_{B(x_{0},r)}(x,z)$ on $B(x_{0},r)$ of
$X^{b}$ is comparable to that of
$X^{\eps_0}$.
Thus for every $x\in B(x_{0},r/2)$, $u(x)$ is comparable to
$u^{*}(x)$. Theorem \ref{HI1} then follows from the
uniform Harnack inequality for mixed stable processes;
see \cite[(3.40)]{CKS3}.        \qed

\begin{lemma}[Harnack inequality]\label{HI2}
Under Assumption 1,
there exists a constant $C_{8}=C_{8}(d,\alpha,\beta, \\ {M_1},{M_2})>0$ such that the following statement is
true: If $x_{1},x_{2}\in \mathbb{R}^{d}$, $r\in (0,r_{1}]$ and
$k\in\mathbb{N}$ are such that $|x_{1}-x_{2}|<2^{k}r$, then for
 every non-negative function $u$ which is harmonic with respect to $X^{b}$
in $B(x_{1},r)\cup B(x_{2},r)$, we have
\begin{equation}
C_{8}^{-1}2^{-k(d+\alpha)}u(x_{2})\le u(x_{1})\le
C_{8}2^{k(d+\alpha)}u(x_{2}).\label{lem10.main}
\end{equation}
\end{lemma}

\proof Without loss of generality, we may assume $|x_{1}-x_{2}|\ge
r/4$. Note that for every $x\in B(x_{2},r/8)\subset
B(x_{1},r/8)^{c}$, we have $|x-x_{1}|<2^{k+1}r$. Thus by Lemma
\ref{lemma2} and Assumption 1, we have
\begin{eqnarray} \label{lem10.1}
K^{b}_{B(x_{1},r/8)}(x_{1},x)&\ge&\frac{1}{2{M_2}}\int_{B(x_{1},r/8)}G_{B(x_{1},r/8)}(x_{1},y)J(y,x)dy\nonumber\\
&=&\frac{1}{2{M_2}}K_{B(x_{1},r/8)}(x_{1},x)\nonumber\\
&\ge&\frac{c_{1}}{2{M_2}}2^{-\alpha}r^{\alpha}(2^{-k-1}r)^{-d-\alpha}=c_{2}r^{-d}2^{-k(d+\alpha)}
.
\end{eqnarray}
Recall that by Theorem \ref{HI1}, we have
$u(x)\ge c_{3}u(x_{2})$ for every
$x\in B(x_{2},r/8)$.
Thus by \eqref{lem10.1},
\begin{eqnarray}
u(x_{1})&\ge&\int_{B(x_{2},r/8)}u(x)K^{b}_{B(x_{1},r/8)}(x_{1},x)dx\nonumber\\
&\ge&c_{2}c_{3}u(x_{2})r^{-d}2^{-k(d+\alpha)}\int_{B(x_{2},r/8)}dx\nonumber\\
&\ge&c_{4}2^{-k(d+\alpha)}u(x_{2}),\nonumber
\end{eqnarray}
and  \eqref{lem10.main} follows by symmetry. \qed

\medskip

\noindent{\bf Proof of Theorem \ref{BHP}.}
Note that there are constants $R_0=R_0 (d, \alpha, \beta, M_1)\in (0,r_{1})$
and $c=c(d, \alpha, \beta, M_1)>1$ so that
\begin{equation}
\frac{c^{-1}}{|x-y|^{d+\alpha}} \leq J^b (x, y) \leq \frac{c}{|x-y|^{d+\alpha}}
 \nonumber
\end{equation}
 for all $|x-y|\leq  R_0$ and
 $b(x, z)$ satisfying \eqref{condi0}.
Thus using \eqref{e:3.3}, we can get
uniform estimates on the Poisson kernel
$$K^{b}_{B(x_{0},r)}(x,z)=\int_{B(x_{0},r)}G^{b}_{B(x_{0},r)}(x,y)J^{b}(y,z)dy$$
of any ball $B(x_{0},r)$ with respect to $X^{b}$ with $r\in (0, R_0/3)$, $x\in B(x_0, r)$
and $r<|z-x_0|<2r$. Specifically, for
$r<|z-x_{0}|<2r$, $K^{b}_{B(x_{0},r)}(x,z)$ is uniformly comparable to
$K_{B(x_{0},r)}(x,z)$. Using the
explicit formula \eqref{ballpoisson} for the Poisson kernel $K_{B(x_0, r)}$,
\eqref{e:3.3},
Theorem \ref{HI1} and \eqref{condi2},
we can adapt the
arguments in \cite[Theorem 2.6]{CKS1} to get our uniform boundary Harnack principle
\ref{BHP} (cf. the proof of \cite[Theorem 3.9]{CKS3}).
Since the proof is almost identical to those in \cite[Section 3]{CKS1},
we omit the details here.           \qed

\begin{lemma}\label{lemma13}
Suppose Assumption 1 holds and $D$ is a Lipschitz open set with characteristics $(\lambda_{0},R_{0})$.
Let $r_1\in (0, 1]$ be the constant
in Lemma \ref{lemma2}.  There is a positive constant
$C_{9}=C_{9}(d,\alpha,\beta,\lambda_{0},R_{0},{M_1},{M_2}) \geq 1$
such that for every $z_{0}\in\partial D$, $r\in (0,r_{1}/2)$, and
every non-negative harmonic function $u$
that is regular harmonic in $D\cap B(z_0, 2r)$ with respect to $X^b$ and
vanishes in $D^c \cap B(z_0, 2r)$,
\begin{equation}\label{lem13.main}
\mathbb{E}_{x}\left[u(X^{b}_{\tau_{D\cap
B_{k}}}):X^{b}_{\tau_{D\cap B_{k}}}\in
B^{c}_{0}\right]\le
C_{9} 2^{-k\alpha} u(x),  \quad    x\in D\cap B_{k},
\end{equation}
for $B_{k}:=B(z_0,2^{-k}r)$ and  $k\geq 1$.
\end{lemma}

\proof Without loss of generality, we may assume $z_0=0$.
By the uniform inner cone property of a Lipschitz open set,
one can find a point $\tilde{z}_{0}\in D\cap B(0,r)$ and
$\kappa=\kappa(\lambda_{0},R_{0}) \in (0,1)$ such that
$\widetilde{B}_{k}:=B(\tilde{z}_{k},\kappa 2^{-k}r)\subset B_{k}\cap D$ for every
$\tilde{z}_{k}:=2^{-k}\tilde{z}_{0}$ and $k\ge 0$. Define
$$u_{k}(x):=\mathbb{E}_{x}\left[u(X^{b}_{\tau_{D\cap B_{k}}}):X^{b}_{\tau_{D\cap B_{k}}}\in
B^{c}_{0}\right].$$
Since $u_{0}=u$, \eqref{lem13.main} is clearly true for $k=0$.
Henceforth we assume $k\ge 1$.
Note that $u_{k}\ge 0$ is
regular harmonic with respect to $X^{b}$ in $D\cap B_{k}$, and
$u_{k}(x)\le u_{k-1}(x)$ for all $x\in \mathbb{R}^{d}$. Define
$$I_{k}(x):=\mathbb{E}_{x}\left[u(X^{b}_{\tau_{B_{k}}}):X^{b}_{\tau_{B_{k}}}\in
B^{c}_{0}\right].$$
Clearly by definition $u_{k}(\tilde{z}_{k})\le I_{k}(\tilde{z}_{k})$. For any $k\ge
1$, by Lemma \ref{lemma2} and \eqref{scalingforgb}, we have
\begin{eqnarray}\label{9}
K^{b}_{B_{k}}(\tilde{z}_{k},y)&=&\int_{B_{k}}G^{b}_{B_{k}}(\tilde{z}_{k},z)J^{b}(z,y)dz\nonumber\\
&\le&\frac{3}{2}\int_{B_{k}}G_{B_{k}}(\tilde{z}_{k},z)J^{b}(z,y)dz\nonumber\\
&=&\frac{3}{2}\int_{2^{-(k-1)}B_{1}}G_{2^{-(k-1)}B_{1}}(2^{-(k-1)}\tilde{z}_{1},z)J^{b}(z,y)dz\nonumber\\
&=&\frac{3}{2}\,2^{-(k-1)d}\int_{B_{1}}G_{2^{-(k-1)}B_{1}}(2^{-(k-1)}\tilde{z}_{1},2^{-(k-1)}w)J^{b}(2^{-(k-1)}w,y)dw\nonumber\\
&=&\frac{3}{2}\,2^{-(k-1)\alpha}\int_{B_{1}}G_{B_{1}}(\tilde{z}_{1},w)J^{b}(2^{-(k-1)}w,y)dw.
\end{eqnarray}
Note that for any
$y\in B^{c}_{0}$
and $w\in B_{1}$,
 $$
 \frac{|y-w|}{|y-2^{-(k-1)}w|}\le \frac{|y|+|w|}{|y|-2^{-k+1}|w|}\le 3 .
 $$
Thus by \eqref{condi2} we have
\begin{eqnarray}
J^{b}(2^{-(k-1)}w,y)\le {M_2}
J^{\eps_0}
(|y-2^{-(k-1)}w|)\le 3^{d+\alpha}{M_2} J^{\eps_0}(|y-w|)
\le 3^{d+\alpha}M^2_2 J^{b}(w,y).
\label{20}
\end{eqnarray}
It follows from \eqref{9}, \eqref{20} and Lemma \ref{lemma2} that
for any $y\in B^{c}_{0}$,
\begin{eqnarray}
K^{b}_{B_{k}}(\tilde{z}_{k},y)&\le&\frac{3^{d+\alpha+1}}{2}M^2_22^{-(k-1)\alpha}\int_{B_{1}}G_{B_{1}}(\tilde{z}_{1},w)J^{b}(w,y)dw\nonumber\\
&\le&3^{d+\alpha+1}M^2_22^{-(k-1)\alpha}\int_{B_{1}}G^{b}_{B_{1}}(\tilde{z}_{1},w)J^{b}(w,y)dw\nonumber\\
&=&c_{1}2^{-k\alpha}K^{b}_{B_{1}}(\tilde{z}_{1},y).
\nonumber
\end{eqnarray}
Now we have for  $k\ge 1$
\begin{eqnarray}
I_{k}(\tilde{z}_{k})&=&\int_{B^{c}_{0}}u(y)K^{b}_{B_{k}}(\tilde{z}_{k},y)dy\nonumber\\
&\le&c_{1}2^{-k\alpha}\int_{B^{c}_{0}}u(y)K^{b}_{B_{1}}(\tilde{z}_{1},y)dy\nonumber\\
&=&c_{1}2^{-k\alpha}I_{1}(\tilde{z}_{1}).\label{3.7}
\end{eqnarray}
Next we   compare $I_{1}(\tilde{z}_{1})$ with $u(\tilde{z}_{1})$.  Using Lemma
\ref{lemma2}, \eqref{condi2} and \eqref{scalingforgb}, we have
\begin{eqnarray}
K^{b}_{\widetilde{B}_{1}}(\tilde{z}_{1},y)&=&\int_{|z-\tilde{z}_{1}|<\kappa
r/2}G^{b}_{\widetilde{B}_{1}}(\tilde{z}_{1},z)J^{b}(z,y)dz\nonumber\\
&\ge&\frac{1}{2{M_2}}\int_{|z-\tilde{z}_{1}|<\kappa
r/2}G_{\widetilde{B}_{1}}(\tilde{z}_{1},z)J^{\eps_0}(|y-z|)dz\nonumber\\
&=&\frac{1}{2{M_2}}\int_{|z-\tilde{z}_{1}|<\kappa
r/2}G_{\kappa B_{1}}(0,z-\tilde{z}_{1})J^{\eps_0}(|y-z|)dz\nonumber\\
&=&\frac{1}{2{M_2}}\int_{|w|<\kappa
r/2}G_{\kappa  B_{1}}(0,w)J^{\eps_0}(|y-\tilde{z}_{1}-w|)dw\nonumber\\
&=&\frac{1}{2{M_2}}\kappa ^{d}\int_{|z|<r/2}G_{\kappa  B_{1}}(0,\kappa  z)J^{\eps_0}(|y-\tilde{z}_{1}-\kappa z|)dz\nonumber\\
&=&\frac{1}{2{M_2}}\kappa^{\alpha}\int_{B_{1}}G_{B_{1}}(0,z)J^{\eps_0}(|y-\tilde{z}_{1}-\kappa
z|)dz.\label{3.8}
\end{eqnarray}
Again using Lemma \ref{lemma2} and \eqref{condi2}, we have
\begin{eqnarray}
K^{b}_{B_{1}}(\tilde{z}_{1},y)&=&\int_{B_{1}}G^{b}_{B_{1}}(\tilde{z}_{1},z)J^{b}(z,y)dz\nonumber\\
&\le&\frac{3}{2}{M_2}\int_{B_{1}}G_{B_{1}}(\tilde{z}_{1},z)J^{\eps_0}(|y-z|)dz\nonumber\\
&=&\frac{3}{2}{M_2}\left(\int_{|z|\le|\tilde{z}_{1}|/2}+\int_{|\tilde{z}_{1}|/2<|z|<r/2}
G_{B_{1}}(\tilde{z}_{1},z)J^{\eps_0}(|y-z|)dz\right)\nonumber\\
&=&\frac{3}{2}{M_2}\left(\int_{|z|\le|\tilde{z}_{1}|/2}G_{B_{1}}
(\tilde{z}_{1},z)J^{\eps_0}(|y-z|)dz\right.\nonumber\\
&&\left.+2^{d}\int_{|\tilde{z}_{1}|/4<|w+\tilde{z}_{1}/2|<r/4}G_{B_{1}}
(\tilde{z}_{1},2w+\tilde{z}_{1})J^{\eps_0}(|y-\tilde{z}_{1}-2w|)dw\right).\label{3.9}
\end{eqnarray}
Note that for any  $y\in B^{c}_{0}$
and $|z|\le |\tilde{z}_{1}|/2$,
$|z-\tilde{z}_{1}|\ge |\tilde{z}_{1}|-|z|\ge |z|$ and\
$|y-\tilde{z}_{1}-\kappa z|/|y-z|\le
\left(|y|+|\tilde{z}_{1}|+\kappa|z|\right)/\left(|y|-|z|\right)\le 4$. Thus
\begin{eqnarray}
G_{B_{1}}(\tilde{z}_{1},z)&\asymp&|z-\tilde{z}_{1}|^{\alpha -d}\left(1\wedge\frac{\delta_{B_{1}}(\tilde{z}_{1})^{\alpha/2}\delta_{B_{1}}(z)^{\alpha/2}}{|z-\tilde{z}_{1}|^{\alpha}}\right)\nonumber\\
&\le&|z|^{\alpha -d}\left(1\wedge\frac{\delta_{B_{1}}(0)^{\alpha/2}\delta_{B_{1}}(z)^{\alpha/2}}
{|z|^{\alpha}}\right)
\asymp G_{B_{1}}(0,z), \nonumber
\end{eqnarray}
and
$$
J^{\eps_0}(|y-z|)\le J^{\eps_0}(\frac{1}{4}|y-\tilde{z}_{1}-\kappa z|)\le
4^{d+\alpha}J^{\eps_0}(|y-\tilde{z}_{1}-\kappa z|).
$$
It follows then that for any  $y\in B^{c}_{0}$,
\begin{equation}
\int_{|z|\le|\tilde{z}_{1}|/2}G_{B_{1}}(\tilde{z}_{1},z)J^{\eps_0}(|y-z|)dz\le
c_{2}\int_{|z|\le|\tilde{z}_{1}|/2}G_{B_{1}}(0,z)J^{\eps_0}(|y-\tilde{z}_{1}-\kappa
z|)dz.\label{3.12}
\end{equation}
Note that for  $y\in B^{c}_{0}$
and $|\tilde{z}_{1}|/4<|w+\tilde{z}_{1}/2|<r/4$,
  $\delta_{B_{1}}(2w+\tilde{z}_{1})=r/2-|2w+\tilde{z}_{1}|\le
2(r/2-|w|)=2\delta_{B_{1}}(w)$, and $|y-\tilde{z}_{1}-\kappa
w|/|y-\tilde{z}_{1}-2w|\le
\left(|y|+\kappa |w+\tilde{z}_{1}/2|+(1-\kappa/2)|\tilde z_1|\right)
/\left(|y|-|\tilde{z}_{1}+2w|\right)
\leq 2
$.
Thus
\begin{eqnarray}
G_{B_{1}}(\tilde{z}_{1},2w+\tilde{z}_{1})&\asymp&|2w|^{\alpha -d}\left(1\wedge\frac{\delta_{B_{1}}(\tilde{z}_{1})^{\alpha/2}\delta_{B_{1}}(2w+\tilde{z}_{1})^{\alpha/2}}{|2w|^{\alpha}}\right)\nonumber\\
&\lesssim&|w|^{\alpha -d}\left(1\wedge\frac{\delta_{B_{1}}(0)^{\alpha/2}\delta_{B_{1}}(w)^{\alpha/2}}
{|w|^{\alpha}}\right)
\asymp G_{B_{1}}(0,w), \nonumber
\end{eqnarray}
and
$$
J^{\eps_0}(|y-\tilde{z}_{1}-2w|)\le
J^{\eps_0}( |y-\tilde{z}_{1}-\kappa w|/2 )\le
2^{d+\alpha}J^{\eps_0}(|y-\tilde{z}_{1}-\kappa w|).
$$
Thus
for any $y\in B_{0}^{c}$,
\begin{eqnarray}
&&\int_{|\tilde{z}_{1}|/4<|w+\tilde{z}_{1}/2|<r/4}G_{B_{1}}(\tilde{z}_{1},2w+\tilde{z}_{1})J^{\eps_0}(|y-\tilde{z}_{1}-2w|)dw\nonumber\\
&\le&c_{3}\int_{|\tilde{z}_{1}|/4<|w+\tilde{z}_{1}/2|<r/4}G_{B_{1}}(0,w)J^{\eps_0}(|y-\tilde{z}_{1}-\kappa w|)dw\nonumber\\
&\le&c_{3}\int_{B_{1}}G_{B_{1}}(0,w)J^{\eps_0}(|y-\tilde{z}_{1}-\kappa
w|)dw.\label{3.15}
\end{eqnarray}
Using \eqref{3.12} and \eqref{3.15}, we can continue the estimates
in \eqref{3.9} to get that for any $y\in B^{c}_{0}$
\begin{equation}
K^{b}_{B_{1}}(\tilde{z}_{1},y)\le
c_{4}\int_{B_{1}}G_{B_{1}}(0,z)J^{\eps_0}(|y-\tilde{z}_{1}-\kappa
z|)dz\label{3.16}
\end{equation}
Combining \eqref{3.8} and \eqref{3.16}, we get
\begin{equation*}
K^{b}_{B_{1}}(\tilde{z}_{1},y)\le c_{5}\kappa
^{-\alpha}K^{b}_{\widetilde{B}_{1}}(\tilde{z}_{1},y),\quad\hbox{for }  y\in B^{c}_{0}.
\end{equation*}
It follows that
\begin{eqnarray}
I_{1}(\tilde{z}_{1})&=&\int_{B^{c}_{0}}u(y)K^{b}_{B_{1}}(\tilde{z}_{1},y)dy
\leq
c_{5}\kappa^{-\alpha}\int_{B^{c}_{0}}u(y)K^{b}_{\widetilde{B}_{1}}(\tilde{z}_{1},y)dy\nonumber\\
&\le&c_{5}\kappa ^{-\alpha}\int_{\widetilde{B}^{c}_{1}}u(y)K^{b}_{\widetilde{B}_{1}}(\tilde{z}_{1},y)dy
 = c_{5}\kappa ^{-\alpha}u(\tilde{z}_{1}).\label{3.18}
\end{eqnarray}
Consequently by \eqref{3.7} and \eqref{3.18} we have for all $k\ge 1$,
\begin{equation}
u_{k}(\tilde{z}_{k})\le I_{k}(\tilde{z}_{k})\le c_{1}c_{5}\kappa
^{-\alpha}2^{-k\alpha}u(\tilde{z}_{1}).\label{3.19}
\end{equation}
By the monotonicity of $u_{k}$ in $k$, Theorem \ref{BHP}, \eqref{3.19}
and Lemma \ref{HI2}, we conclude that for any $x\in D\cap B_{k}$ and
$k\ge 1$
\begin{equation}
\frac{u_{k}(x)}{u(x)}\le \frac{u_{k-1}(x)}{u(x)}\le
c_{6}\frac{u_{k-1}(\tilde{z}_{k-1})}{u(\tilde{z}_{k-1})}\le c_{6}c_{1}c_{5}\kappa
^{-\alpha}2^{-(k-1)\alpha}\frac{u(\tilde{z}_{1})}{u(\tilde{z}_{k-1})}\le
c_{7}2^{-k\alpha}.\nonumber
\end{equation}
The proof is now complete. \qed

\medskip

For a Lipschitz open set $D$ with characteristics $(\lambda_{0},R_{0})$,
let $\kappa=\kappa (\lambda_0, R_0)\in (0, 1)$ so that $D$ is $\kappa$-fat.
For $z_{0}\in \partial D$ and $r\in (0,1]$,  we use $A_{r}(z_{0})$ to denote
 a point in $D$ such that $B(A_{r}(z_{0}),\kappa r)\subset D\cap B(z_{0},r)$.

 The following lemma follows from Theorem \ref{BHP} and Lemma
\ref{lemma13}
(in place of
 \cite[Lemma 13 and Lemma 14]{Bogdan1})
in the same way as for the case of symmetric
$\alpha$-stable process in \cite[Lemma 16]{Bogdan1}.
We omit its proof here. See also \cite[Lemma 3.4]{KSV} for detailed computation
of a similar result for symmetric L\'evy processes.

\begin{lemma}\label{lemma14}
Suppose Assumption 1 holds and $D$ is a Lipschitz open set with characteristics $(\lambda_{0},R_{0})$.
Let $r_1\in (0, 1]$ be the constant in Lemma \ref{lemma2}.
There exist positive constants
${\gamma_1}={\gamma_1}(d,\alpha,\beta,\lambda_{0},R_{0},{M_1},{M_2}) $ and
$C_{10}=C_{10}(d,\alpha,\beta,\lambda_{0},R_{0},{M_1}, {M_2}) $  such that
for every $z_{0}\in \partial D$, $r\in (0,r_{1}/2)$ and all non-negative functions $u,\ v$ that are regular harmonic in $D\cap B(z_{0},2r)$ and vanish in $D^{c}\cap B(z_{0},2r)$ with $u(A_{r}(z_{0}))=v(A_{r}(z_{0}))>0$, we have
\begin{description}
\item{\rm (i)} $h(z_0):=\lim_{D \ni x\to z_0}u(x)/v(x)$ exists;

\item{\rm (ii)} $|\frac{u(x)}{v(x)}-h(z_0)|\le C_{10}\left(\frac{|x-z_0|}{r}\right)^{{\gamma_1}}   \hbox{ for }  x\in D\cap B(z_0,r)$.
\end{description}
\end{lemma}

\section{Gradient upper bound estimates}\label{S:4}

We now study gradient estimates for non-negative harmonic functions of $X^b$
in open sets.

\begin{lemma}\label{lemma4}
Suppose $b$ is a bounded function satisfying \eqref{condi0},
  and that for every $x\in\mathbb{R}^{d}$,
$J^{b}(x,y)\ge 0$ a.e. $y\in \mathbb{R}^{d}$.
Let $r_1 \in (0, 1]$ be the constant in Lemma \ref{lemma2},
and  $B=B(x_0, r)$ with $r\in (0, r_1]$. Then
for every $x\in B$, $z\in \bar{B}^{c}$ and $1\leq i\leq d$,
\begin{equation}
\partial_{x_i}\int_{B}G_{B}(x,y)J^{b
}(y,z)dy=\int_{B}\partial_{x_i}G_{B}(x,y)J^{b
}(y,z)dy,\label{lem9.1}
\end{equation}
\begin{equation}
\partial_{x_i}\int_{B}\left(\int_{B}G_{B}(x,y)\mathcal{S}^{b
}_{y}G^{b }_{B}(y,w)dy\right) J^{b
}(w,z)dw=\int_{B}\left( \int_{B}\partial_{x_i}G_{B}(x,y)\mathcal{S}^{b
}_{y}G^{b }_{B}(y,w)dy\right) J^{b }(w,z)dw, \label{lem9.2}
\end{equation}
and
\begin{equation} \label{2.21}
\partial_{x_i}K^{b}_{B}(x,z)
=\int_{B}\partial_{x_i}G_{B}(x ,y)J^{b}(y,z)dy
+\int_{B}\left(\int_{B}\partial_{x_i}G_{B}(x
,y)\mathcal{S}^{b}_{y}G^{b}_{B}(y,w)dy\right) J^{b}(w,z)dw.
\end{equation}
\end{lemma}

\proof Without loss of generality we assume $i=d$. Fix $x\in B$
and $z\in \bar{B}^{c}$.   We have
$$\sup_{y\in B}|J^b(y, z)|\le
\frac{\mathcal{A}(d,-\alpha)}{\delta_{B}(z)^{d+\alpha}}+\frac{\|b
\|_{\infty}\mathcal{A}(d,-\beta)}{\delta_{B}(z)^{d+\beta}}<+\infty.$$
Thus \eqref{lem9.1} follows directly from
\cite[Lemma 5.2]{BKN}.

Let $g_{z}(y):=\int_{B}\mathcal{S}^{b }_{y}G^{b}_{B}(y,w)J^{b }(w,z)dw$
for $y\in B$.  We have
\begin{equation}
\partial_{x_d}\int_{B}G_{B}(x,y)g_{z}(y)dy=\lim_{\lambda\to
0}\int_{B}\left[\frac{G_{B}(x+\lambda e_{d},y)-G_{B}(x,y)}{\lambda}\right]g_{z}(y)dy.\label{2.13}
\end{equation}
To prove \eqref{lem9.2}, we only need to show that the integrand in
the right hand side of \eqref{2.13} is uniformly  integrable
on $B$ in $\lambda\in (0, \delta_{B}(x)/2)$. Note that we have
$$
G_{B}(x,y)=G(x,y)-\mathbb{E}_{y}\left[G(x, X_{\tau_{B}})\right]=:G(x, y)-H(x,y) .
$$
Thus
\begin{eqnarray}
\frac{|G_{B}(x+\lambda e_d,y)-G_{B}(x,y)|}{\lambda}&\le&\frac{|G(x+\lambda e_d,y)-G(x,y)|}{\lambda}+\frac{|H(x+\lambda e_d,y)-H(x,y)|}{\lambda}\nonumber\\
&=:&I+II.\nonumber
\end{eqnarray}
Obviously by \eqref{globalgreen} we have
\begin{equation}
I\le c_{1}\left(|x+\lambda e_{d}-y|^{\alpha -1 -d}+
|x-y|^{\alpha-d-1} \right) \quad\hbox{for }
 y\in B,
 \nonumber
\end{equation}
for some positive constant $c_{1}=c_{1}(d,\alpha)$.
Since $H(x, y)= \E_y [ G(x, X_{\tau_B)}]$,
by the mean-value theorem, there is a point $x_\lambda$ in the line segment connecting
$x$ with $x+{\lambda} e_d$ so that
$$
II= \partial_{x_d} H(x_\lambda, y)=  \E_y  [ \partial_{x_d} G(x_\lambda, X_{\tau_B})] \leq c_{2} \delta_B (x)^{\alpha -1 -d}.
$$
Thus for some positive constant $c_{3}=c_{3}(d,\alpha,x)$, we have
\begin{equation}
\frac{|G_{B}(x+\lambda e_d,y)-G_{B}(x,y)|}{\lambda}\le
c_{1}\left(|x+\lambda e_d-y|\wedge
|x-y|\right)^{\alpha-d-1}+c_{3}.\label{2.16}
\end{equation}
Let $h(y):=\int_{B}h_{B}(y,w)dw$ for $y\in B$.
Note that by Lemma \ref{lemma2} and the boundedness of
$w\mapsto J^b (w, z) $ on $B$,
\begin{equation}\label{2.17}
|g_{z}(y)|\le c_{4}\int_{B}h_{B}(y,w) J^b (w, z)  dw\le c_{5}h(y).
\end{equation}
Thus by \eqref{2.16} and \eqref{2.17} the integrand in the right
hand side of \eqref{2.13} is uniformly  integrable on $B$ in
$\lambda\in (0, \delta_B(x)/2)$
if the following three conditions are true:
\begin{description}
\item{(i)} $\int_{B}h(y)dy<+\infty$;

\item{(ii)} $\sup_{w\in
B(x,\delta_{B}(x)/2)}\int_{B}h(y)|y-w|^{\alpha -1-d}\,dy<+\infty$;

\item{(iii)} $\lim_{\varepsilon\downarrow 0}\sup_{w\in B(x,\delta_{B}(x)/2)}\int_{\{y\in
B:|y-w|<\varepsilon\}}h(y)|y-w|^{\alpha -1-d}=0$.
\end{description}
If $\alpha>2 \beta$, then for any $y\in B$,
$$h(y)\le
\int_{w\in B}|y-w|^{\alpha-\beta-d}dw\le
\int_{|u|<2r}|u|^{\alpha-\beta-d}du<+\infty;
$$
that is, $h(y)$ is bounded from above on $B$.  Obviously (i)-(iii) hold for $h$.
If $\alpha=2\beta$, then
\begin{eqnarray}
h(y)&=&\int_{w\in
B}|w-y|^{\beta-d}\left(1\wedge \frac{\delta_{B}(w)^{\beta}}{|y-w|^{\beta}}\right) \left(1\vee
\log\frac{|w-y|}{\delta_{B}(y)}\right)dw\nonumber\\
&=&\int_{w\in
B,|w-y|>e\delta_{B}(y)}|w-y|^{\beta/2-d}\frac{\delta_{B}(w)^{\beta}}{\delta_{B}(y)^{\beta/2}}\left(\frac{\delta_{B}(y)^{\beta/2}}{|w-y|^{\beta/2}}
\log\frac{|w-y|}{\delta_{B}(y)}\right)dw\nonumber\\
&&+\int_{w\in B,|w-y|\le
e\delta_{B}(y)}|w-y|^{\beta -d } dw\nonumber\\
&\lesssim&\delta_{B}(y)^{-\beta/2}+1.\label{2.18}
\end{eqnarray}
Using this upper bound, it is easy to check $h$ satisfies (i) and (ii).
 As for (iii), note that
 $\delta_{B}(w)\ge \delta_{B}(x)/2$ for every $w\in B(x,\delta_{B}(x)/2)$.
 Consider an arbitrary $\varepsilon\in
(0,\delta_{B}(x)/4)$. Then $B(w,\varepsilon)\subset B$  and
$\delta_{B}(y)\ge \delta_{B}(x)/4$ for every $y\in
B(w,\varepsilon)$.
We have by \eqref{2.18},
\begin{eqnarray}
\int_{|y-w|<\varepsilon}h(y)|y-w|^{\alpha-1-d}dy&\lesssim&\int_{|y-w|<\varepsilon}
(\delta_{B}(y)^{-\beta/2}+1)|w-y|^{\alpha-1-d}dy\nonumber\\
&\lesssim&\int_{|y-w|<\varepsilon}(\delta_{B}(x)^{-\beta/2}+1)|w-y|^{\alpha-1-d}dy.\nonumber
\end{eqnarray}
Thus condition (iii)  is implied by
the fact that
$$
\lim_{\varepsilon\downarrow 0}\sup_{w\in B(x,\delta_{B}(x)/2)}
\int_{|y-w|<\varepsilon}(\delta_{B}(x)^{-\beta/2}+1)|w-y|^{\alpha-1-d}dy=0.
$$
When  $\alpha <2\beta$,  similar to \eqref{2.18}
we have
$$
h(y)\asymp\int_{w\in B}|w-y|^{\beta-d}
\left( 1 \wedge \frac{\delta_{B}(w)^{\beta}}{|y-w|^{\beta}}\right)
\left(1\vee\frac{|y-w|^{\beta-\alpha/2}}{\delta_{B}(y)^{\beta-\alpha/2}}\right)dw\lesssim
\delta_{B}(y)^{\beta-\alpha/2}+1.$$
By a similar
calculations as in the case $\alpha=2 \beta$,
we can show that  (i)-(iii) hold
for $h$.  This completes the proof.
\qed

\bigskip

\begin{lemma}\label{lemma6}
Under Assumption 1,
there exists a constant $C_{11}=C_{11}(d,\alpha,\beta, {M_1}, {M_2})>0$
such that for every $B=B(x_0, 1)$ and $1\leq i\leq d$,
\begin{equation}
\int_{B}|\partial_{x_i}G_{B}(x, y)|J^{b}(y,z)dy\le
C_{11}\int_{B}G_{B}(x_0,y)J^{b}(y,z)dy\quad\hbox{for }
x\in B(x_{0},1/4)\mbox{ and } z\in\bar{B}^{c}.\label{lem6.main}
\end{equation}
\end{lemma}

\proof Without loss of generality, we  assume $x_0=0$ and $i=d$.
For every $|x|<1/4$ and $|y|<1$, we have $|x-y|\wedge
\delta_{B}(x)\asymp |x-y|$. Thus by \eqref{partialgreen},
\begin{eqnarray}
\int_{B}|\partial_{x_d}G_{B}(x,y)|J^{b}(y,z)dy
&\lesssim&\int_{|y|<1}\frac{G_{B}(x,y)}{|x-y|\wedge
\delta_{B}(x)}\,J^{b}(y,z)dy\nonumber\\
&\asymp&\int_{|y|<1}\frac{G_{B}(x,y)}{|x-y|}\,J^{b}(y,z)dy\nonumber\\
&=& \left( \int_{1/2\le |y|<1}+\int_{|y|<1/2} \right) \frac{G_{B}(x,y)}{|x-y|}\,J^{b}(y,z)dy\nonumber\\
&=:&I(x,z)+II(x,z).
\nonumber
\end{eqnarray}
For $1/2\le |y|<1$ and $|x|<1/4$, we have $|x-y|\asymp |y|\asymp
1$ and $\delta_{B}(x)\asymp 1$. Thus
$$G_{B}(x,y)\asymp |y|^{\alpha -d}\left(1\wedge\frac{\delta_{B}(y)^{\alpha/2}}{|y|^{\alpha}}\right)\asymp G_{B}(0,y),$$
and consequently
\begin{eqnarray}\label{2.27}
I(x,z) \asymp\int_{1/2\le|y|<1}G_{B}(0,y)J^{b}(y,z)dy.
\end{eqnarray}
For every $|y|<1/2$, $|x|<1/4$ and $|z|>1$, we have
$|z-y|>1/2$, $|z-y|\asymp |z-y+x|$ and
$\delta_{B}(y)=1-|y|\asymp 1-|y-x|\asymp 1$.
Thus by \eqref{condi2}
\begin{eqnarray}
II(x,z)&\stackrel{c_{1}({M_2})}{\asymp}&\int_{|y|<1/2}|x-y|^{-d+\alpha-1}\left(1\wedge
\frac{\delta_{B}(x)^{\alpha/2}\delta_{B}(y)^{\alpha/2}}{|x-y|^{\alpha}}\right)
J^{\eps_0}(|z-y|)dy\nonumber\\
&\asymp&\int_{|y|<1/2}|x-y|^{-d+\alpha-1}\left(1\wedge
\frac{\delta_{B}(x)^{\alpha/2}(1-|y-x|)^{\alpha/2}}{|x-y|^{\alpha}}\right)
J^{\eps_0}(|z-y+x|)dy\nonumber\\
&\asymp&\int_{|w+x|<1/2}|w|^{-d+\alpha-1}\left(1\wedge
\frac{(1-|w|)^{\alpha/2}}{|w|^{\alpha}}\right)J^{\eps_0}(|z-w|)dw\nonumber\\
&\le&\int_{|w|<3/4}|w|^{-d+\alpha-1}\left(1\wedge
\frac{(1-|w|)^{\alpha/2}}{|w|^{\alpha}}\right)J^{\eps_0}(|z-w|)dw\nonumber\\
&\stackrel{c_{2}({M_2})}{\asymp}&\int_{|w|<3/4}|w|^{-d+\alpha-1}\left(1\wedge
\frac{(1-|w|)^{\alpha/2}}{|w|^{\alpha}}\right)J^{b}(w,z)dw\nonumber\\
&=:&g_{1}(z).
\nonumber
\end{eqnarray}
For every $z>1$, let
$$g_{2}(z):=\int_{|w|<3/4}|w|^{\alpha -d}\left(1\wedge
\frac{(1-|w|)^{\alpha/2}}{|w|^{\alpha}}\right)J^{b}(w,z)dw.$$
Obviously
\begin{equation}
g_{2}(z)\asymp\int_{|w|<3/4}G_{B}(0,w)J^{b}(w,z)dw\le\int_{B}G_{B}(0,w)J^{b}(w,z)dw.\label{2.29}
\end{equation}
Note that $J^{\eps_0}(|x-y|)$ is non-increasing in
$|x-y|$. Thus by \eqref{condi2}
\begin{equation}
\sup_{|z|>1}\frac{g_{1}(z)}{g_{2}(z)}\le\sup_{|z|>1}\frac{{M_2}
J^{\eps_0}(|z|-3/4)\int_{|w|<3/4}|w|^{\alpha-d-1}\left(1\wedge
\frac{(1-|w|)^{\alpha/2}}{|w|^{\alpha}}\right)dw}{M^{-1}_{2}
J^{\eps_0}(|z|+3/4)\int_{|w|<3/4}|w|^{\alpha -d}\left(1\wedge
\frac{(1-|w|)^{\alpha/2}}{|w|^{\alpha}}\right)dw}\le
M<+\infty,\label{2.30}
\end{equation}
where $M=M(d,\alpha,\beta,{M_2})>0$. Thus by \eqref{2.29} and
\eqref{2.30} we prove that
\begin{equation}
II(x,z)\stackrel{c_{3}({M_2})}{\lesssim}\int_{B}G_{B}(0,w)J^{b}(z,w)dw\quad\hbox{for } |x|<1/4,\
|z|>1.\label{2.31}
\end{equation}
Therefore \eqref{lem6.main} follows from \eqref{2.27} and
\eqref{2.31}. \qed

\medskip

Recall the definition of $r_D(x, y)$ and $h_D(x, y)$
from \eqref{e:1.3} and \eqref{e:3.2}, respectively.

\begin{lemma}\label{lemma7}
For $B=B(0,1)$, there exists a constant $C_{12}=C_{12}(d,\alpha,\beta)>0$ such that
for every $1\leq i\leq d$,
$|x|<1/4$ and $|w|<1$,
\begin{equation}
\int_{B}|\partial_{x_i}G_{B}(x,y)|h_{B}(y,w)dy\le C_{12}
|x-w|^{-d+(\alpha-1)\wedge(\alpha-\beta)}\delta_{B}(w)^{\alpha/2} .
\nonumber
\end{equation}
\end{lemma}

\proof Without loss of generality, we assume $i=d$.
For any $|x|<1/4$ and $|y|<1$, we have
$\delta_{B}(x)\asymp 1$ and $|x-y|\wedge
\delta_{B}(x)\asymp |x-y|$. Thus by \eqref{partialgreen},
\begin{eqnarray}
\int_{B}|\partial_{x_d}G_{B}(x,y)|h_{B}(y,w)dy&\lesssim&\int_{|y|<1}\frac{G_{B}(x,y)}{|x-y|\wedge
\delta_{B}(x)}\,h_{B}(y,w)dy\nonumber\\
&\asymp&\int_{|y|<1}\frac{G_{B}(x,y)}{|x-y|}\,h_{B}(y,w)dy.\label{eq1}
\end{eqnarray}
We note that for $|x|<1/4$ and $|y|<1$,
\begin{equation}
G_{B}(x,y)\asymp |x-y|^{\alpha -d}\left(1\wedge
\frac{\delta_{B}(x)^{\alpha/2}}{|x-y|^{\alpha/2}}\right)\left(1\wedge
\frac{\delta_{B}(y)^{\alpha/2}}{|x-y|^{\alpha/2}}\right)
\asymp
|x-y|^{\alpha -d}\frac{\delta_{B}(y)^{\alpha/2}}{r_{B}(x,y)^{\alpha}},\label{2.34}
\end{equation}
and $r_{B}(x,y)\ge\delta_{B}(x)\ge 3/4$.
Now we
calculate the integral in \eqref{eq1} using \eqref{2.34} and the
explicit formula of $h_{B}(y,w)$. If $\alpha>2 \beta$, we have
\begin{eqnarray}
 \eqref{eq1}
&\asymp&\int_{|y|<1}\frac{\delta_{B}(w)^{\alpha/2}}{|x-y|^{d-\alpha+1}|w-y|^{d-\alpha+\beta}}\,\frac{\delta_{B}(y)^{\alpha/2}}{r_{B}(x,y)^{\alpha}r_{B}(y,w)^{\alpha/2}} \,dy\nonumber\\
&\lesssim&\delta_{B}(w)^{\alpha/2}\int_{|y|<1}|x-y|^{-d+\alpha-1}|y-w|^{-d+\alpha-\beta}dy\nonumber\\
&\lesssim&|x-w|^{-d+(\alpha-1)\wedge(\alpha-\beta)}\delta_{B}(w)^{\alpha/2}\nonumber\\
&\le&|x-w|^{-d+(\alpha-1)\wedge(\alpha-\beta)}\delta_{B}(w)^{\beta}.\label{>}
\end{eqnarray}
If $\alpha = 2\beta$, we have by   \eqref{e:2.11} and \eqref{e:3.2}
\begin{eqnarray}
&&\eqref{eq1}\nonumber\\
&\asymp&\int_{|y|<1}|x-y|^{2\beta-1-d}\frac{\delta_{B}(y)^{\beta}}{r_{B}(x,y)^{2\beta}}\,
|y-w|^{\beta-d}
\frac{\delta_{B}(w)^{\beta}}{r_{B}(y,w)^{\beta}}\left(1\vee
\log\frac{|y-w|}{\delta_{B}(y)}\right)dy\nonumber\\
&=&\int_{|y|<1,|y-w|\le
e\delta_{B}(y)}\frac{\delta_{B}(w)^{\beta}}{|x-y|^{d+1-2\beta}|y-w|^{d-\beta}}\,\frac{\delta_{B}(y)^{\beta}}{r_{B}(x,y)^{2\beta}r_{B}(y,w)^{\beta}}  \,dy\nonumber\\
&&+\int_{|y|<1,|y-w|>
e\delta_{B}(y)}\frac{\delta_{B}(w)^{\beta}}{|x-y|^{d+1-2\beta}|w-y|^{d-\beta}}\,
\frac{|y-w|^{\beta/2}\delta_{B}(y)^{\beta/2}}
{r_{B}(x,y)^{2\beta}r_{B}(y,w)^{\beta}}\left(\frac{\delta_{B}(y)^{\beta/2}}{|y-w|^{\beta/2}}
\log\frac{|y-w|}{\delta_{B}(y)}\right)dy\nonumber\\
&\lesssim&\delta_{B}(w)^{\beta}\int_{|y|<1,|y-w|\le
e\delta_{B}(y)}|x-y|^{2\beta-1-d}|y-w|^{\beta-d}
dy\nonumber\\
&&+\delta_{B}(w)^{\beta}\int_{|y|<1,|y-w|>
e\delta_{B}(y)}|x-y|^{2\beta-1-d}|y-w|^{2\beta-d}
dy\nonumber\\
&\lesssim&\delta_{B}(w)^{\beta}|x-w|^{-d+2\beta-1}=|x-w|^{-d+(\alpha-1)\wedge
(\alpha-\beta)}\delta_{B}(w)^{\beta}.\label{=}
\end{eqnarray}

If $\alpha<2\beta$, we have
\begin{eqnarray}
\eqref{eq1}
&\asymp&\int_{|y|<1}|x-y|^{-d+\alpha-1}\frac{\delta_{B}(y)^{\alpha/2}}{r_{B}(x,y)^{\alpha}}\,|y-w|^{-d+\alpha-\beta}
\frac{\delta_{B}(w)^{\alpha/2}}{r_{B}(y,w)^{\alpha/2}}\left(1\vee
\frac{|y-w|^{\beta-\alpha/2}}{\delta_{B}(y)^{\beta-\alpha/2}}\right)dy\nonumber\\
&\le&\int_{|y|<1,|y-w|\ge
\delta_{B}(y)}\frac{\delta_{B}(w)^{\alpha/2}}{|x-y|^{d-\alpha+1}|w-y|^{d-\alpha+\beta}}\,\frac{|w-y|^{\beta-\alpha/2}\delta_{B}(y)^{\alpha-\beta}}
{r_{B}(y,w)^{\alpha/2}}\,dy\nonumber\\
&&+\int_{|y|<1,|y-w|<
\delta_{B}(y)}\frac{\delta_{B}(w)^{\alpha/2}}{|x-y|^{d-\alpha+1}|w-y|^{d-\alpha+\beta}}\,\frac{\delta_{B}(y)^{\alpha/2}}{r_{B}(x,y)^{\alpha}r_{B}(y,w)^{\alpha/2}}\,dy\nonumber\\
&\lesssim&\delta_{B}(w)^{\alpha/2}\int_{|y|<1,|y-w|\ge
\delta_{B}(y)}|x-y|^{-d+\alpha-1}|y-w|^{-d+\alpha-\beta}dy\nonumber\\
&&+\delta_{B}(w)^{\alpha/2}\int_{|y|<1,|y-w|<
\delta_{B}(y)}|x-y|^{-d+\alpha-1}|y-w|^{-d+\alpha-\beta}dy\nonumber\\
&=&\delta_{B}(w)^{\alpha/2}\int_{|y|<1}|x-y|^{-d+\alpha-1}|y-w|^{-d+\alpha-\beta}dy\nonumber\\
&\lesssim&\delta_{B}(w)^{\alpha/2}|x-w|^{-d+(\alpha-1)\wedge(\alpha-\beta)}\int_{|y|<1}(|x-y|^{\alpha-1-(\alpha-1)\wedge(\alpha-\beta)}+|y-w|^{\alpha-\beta-(\alpha-1)\wedge(\alpha-\beta)})dy\nonumber\\
&\lesssim&|x-w|^{-d+(\alpha-1)\wedge(\alpha-\beta)}\delta_{B}(w)^{\alpha/2}.\label{<}
\end{eqnarray}
Lemma \ref{lemma7} follows from \eqref{=}, \eqref{<} and \eqref{>}.
\qed

\begin{lemma} \label{lemma8}
Under Assumption 1,
there exists a constant
$C_{13}=C_{13}(d,\alpha,\beta,{M_1},M_{2})>0$ such that
for $B=B(x_0, 1)$,  $1\leq i\leq d$,
$x\in B(x_0,1/4)$ and $z\in \bar{B}^{c}$,
\begin{equation}
\int_{B}\left[\int_{B}|\partial_{x_i}G_{B}(x,y)\mathcal{S}^{b}_{y}G^{b}_{B}(y,w)|dy\right]J^{b}(w,z)dw\le
C_{13} \int_{B}G_{B}(x_0 ,w)J^{b}(w,z)dw.\label{2.38}
\end{equation}
\end{lemma}

\proof Without loss of generality, we assume that  $x_0=0$ and $i=d$.
Let $r_1\in (0, 1]$ be the constant in Lemma \ref{lemma2}.
By Lemma \ref{lemma2} and the scaling property, we have for $y,w\in B$,
$$|S^{b}_{y}G^{b}_{B}(y,w)|=r_{1}^{d}|S^{b_{r_{1}}}_{y}G^{b_{r_{1}}}_{r_{1}B}(r_{1}y,r_{1}w)|\le c_{1}
r_{1}^{d}|h_{r_{1}B}(r_{1}y, r_{1}w)|=c_{1}r_{1}^{\alpha-\beta} h_{B}(y,w)\le c_{1}h_{B}(y,w).$$
Here $c_{1}=c_{1}(d,\alpha,\beta,M_{1})>0$.
Hence to prove \eqref{2.38}, it suffices to prove that
for $x\in B(0,1/4)$ and  $z\in\bar{B}^{c}$,
\begin{equation}\label{2.39}
\int_{B}\left[\int_{B}|\partial_{x_d}G_{B}(x,y)|h_{B}(y,w)dy\right]J^{b}(w,z)dw\le
c_{2}\int_{B}G_{B}(x_0,w)J^{b}(w,z)dw
\end{equation}
for some $c_{2}=c_{2}(d,\alpha,\beta,{M_1},{M_2})>0$. By Lemma \ref{lemma7}, we have
\begin{eqnarray}
&&\int_{B}\left[\int_{B}|\partial_{x_d}G_{B}(x,y)|h_{B}(y,w)dy\right]J^{b}(w,z)dw\nonumber\\
&\lesssim&\int_{|w|<1}\delta_{B}(w)^{\alpha/2}|w-x|^{-d+(\alpha-1)\wedge(\alpha-\beta)}J^{b}(w,z)dw\nonumber\\
&=&\int_{|w|\le 1/2}+\int_{1/2<|w|<1}\delta_{B}(w)^{\alpha/2}|w-x|^{-d+(\alpha-1)\wedge(\alpha-\beta)}J^{b}(w,z)dw\nonumber\\
&=:&I(x,z)+II(x,z).\nonumber
\end{eqnarray}
Fix $|x|<1/4$ and $|z|>1$. For $1/2<|w|<1$, we have
$|w-x|\asymp |w|\asymp 1$, and consequently $G_{B}(0,w)\asymp
\delta_{B}(w)^{\alpha/2}$. Thus
\begin{equation}
II(x,z)\asymp \int_{1/2<|w|<1}\delta_{B}(w)^{\alpha/2}J^{b}(w,z)dw\asymp\int_{1/2<|w|<1}G_{B}(0,w)J^{b}(w,z)dw.\label{lem8.II}
\end{equation}
For any $|w|\le 1/2$, we have $1/2\le\delta_{B}(w)\le 1$,
$|z-w|\ge 1/2$ and $|z-w+x|\asymp |z-w|$. Hence by \eqref{condi2}
\begin{eqnarray}
I(x,z)&\stackrel{{M_2}}{\asymp}&\int_{|w|\le
1/2}|w-x|^{-d+(\alpha-1)\wedge(\alpha-\beta)}J^{\eps_0}(|z-w|)dw\nonumber\\
&\asymp&\int_{|w|\le
1/2}|w-x|^{-d+(\alpha-1)\wedge(\alpha-\beta)}J^{\eps_0}(|z-w+x|)dw\nonumber\\
&=&\int_{|v+x|\le
1/2}|v|^{-d+(\alpha-1)\wedge(\alpha-\beta)}J^{\eps_0}(|z-v|)dv\nonumber\\
&\le&\int_{|v|\le
3/4}|v|^{-d+(\alpha-1)\wedge(\alpha-\beta)}J^{\eps_0}(|z-v|)dv
 =:
g_{1}(z).   \label{lem8.3}
\end{eqnarray}
We first consider the case $|z|>2$. Let
$g_{2}(z):=\int_{|v|\le3/4}|v|^{\alpha -d}J^{\eps_0}(|z-v|)dv$.
Note that for any $|v|\le 3/4$, we have $G_{B}(0,v)\asymp
|v|^{\alpha -d}$. Thus
\begin{equation}
g_{2}(z)\asymp \int_{|v|\le
3/4}G_{B}(0,v)J^{\eps_0}(|z-v|)dv\le {M_2}
\int_{|v|\le
3/4}G_{B}(0,v)J^{b}(v,z)dv.\label{lem8.1}
\end{equation}
In addition since $J^{\eps_0}(|y|)$ is
non-increasing in $|y|$, we have
\begin{equation} \label{lem8.2}
\sup_{|z|>2}\frac{g_{1}(z)}{g_{2}(z)}\le
\sup_{|z|>2}\frac{J^{\eps_0}(|z|-3/4)\int_{|v|\le
3/4}|v|^{-d+(\alpha-1)\wedge(\alpha-\beta)}dv}
{J^{\eps_0}(|z|+3/4)\int_{|v|\le
3/4}|v|^{\alpha -d}dv}\le M<+\infty ,
\end{equation}
where $M=M(d,\alpha,\beta)>0$. Therefore by \eqref{lem8.3},
\eqref{lem8.1} and \eqref{lem8.2} we have
\begin{equation}
I(x,z)\stackrel{c_{3}({M_2})}{\lesssim}\int_{B}G_{B}(0,w)J^{b}(w,z)dw\quad\hbox{for }  |x|<1/4\mbox{ and } |z|>2.\label{eq2}
\end{equation}
On the other hand if $1<|z|\le 2$, we have
$0<\delta_{B}(z)\le 1$, and by \eqref{ballpoisson}
\begin{equation}
\int_{B}G_{B}(0,w)J^{b}(w,z)dw\ge
M^{-1}_2\int_{B}G_{B}(0,w)J(w,z)dw
= M^{-1}_2
K_{B}(0,z)\asymp
M^{-1}_2\delta_{B}(z)^{-\alpha/2}\ge
M^{-1}_2.\label{lem8.4}
\end{equation}
Note that $|z-w|\ge 1/4$ for any $|w|\le 3/4$. Thus
\begin{equation}
g_{1}(z)\le J^{\eps_0}(1/4)\int_{|w|\le
3/4}|w|^{-d+(\alpha-1)\wedge(\alpha-\beta)}dw
\lesssim 1.\label{lem8.5}
\end{equation}
Thus by \eqref{lem8.3}, \eqref{lem8.4} and \eqref{lem8.5} we have
\begin{equation}
I(x,z)\stackrel{c_{4}({M_2})}{\lesssim}\int_{B}G_{B}(0,w)J^{b}(w,z)dw\quad\hbox{for }
|x|<1/4\mbox{ and } 1<|z|\le 2.\label{lem8.I}
\end{equation}
Now \eqref{2.39} follows from \eqref{lem8.I} and
\eqref{lem8.II}.\qed

\begin{thm}\label{lemma5}
Let $r_1\in (0, 1]$ be the constant in Lemma \ref{lemma2}.
Under Assumption 1,
there exists a constant $C_{14}=C_{14}(d,\alpha,\beta,{M_1},{M_2})>0$
such that for every ball $B_{r}=B(x_0,r)$ with radius $r\in
(0,r_{1}]$ and $1\leq i\leq d$,
\begin{equation}\label{partialk}
|\partial_{x_i}K^{b}_{B_r}(x ,z)|\le
\frac{C_{14}}{r}\,K^{b}_{B_{r}}(x_0,z),\quad\hbox{for } x\in
B(x_0,r/4) \hbox{ and }  z\in \bar{B_r}^{c}.
\end{equation}
\end{thm}

\proof Let $\lambda:=1/r\geq 1/r_1\geq 1$ and define
$b_\lambda (x, z) = \lambda^{\beta -\alpha} b(\lambda^{-1} x, \lambda^{-1} z)$.
Observe that $\| b_\lambda\|_\infty = r^{\alpha -\beta} \|b\|_\infty
\leq r_1^{\alpha -\beta} {M_1}\leq {M_1}$.
By the scaling properties \eqref{e:2.6a} and \eqref{scallingpoisson},
$b_\lambda (x,z)$ satisfies Assumption 1 and
it suffices to show that for
the ball $B=B(x_{0},1)$,
\begin{equation} \label{partialK}
|\partial_{x_i}K^{b_\lambda}_{B}(x, z)|\le C_{14}
K^{b_\lambda}_{B}(x_0, z)  \qquad\hbox{for }  x\in B(x_0 ,1/4) \hbox{ and }
z\in \bar{B}^{c}.
\end{equation}
We know from  \cite[Lemma 4.10]{CY} that
$$G^{b_{\lambda}}_{B}(x,y)=G_{B}(x,y)+\int_{B}G_{B}(x,z)\mathcal{S}^{b_{\lambda}}_{z}G^{b_{\lambda}}_{B}(z,y)dz\quad \hbox{for }  x,y\in B.$$
Thus by \eqref{(3)}, for $i=1,\cdots,d$, every $x\in B$, and
$z\in \bar{B}^{c}$,
\begin{equation}
\partial_{x_i}K^{b_{\lambda}}_{B}(x,z)=\partial_{x_i}\int_{B}G_{B}(x,y)J^{b_{\lambda}}(y,z)dy
+\partial_{x_i}\int_{B}\left(\int_{B}G_{B}(x,y)\mathcal{S}^{b_{\lambda}}_{z}G^{b_{\lambda}}_{B}(y,w)dy\right)
J^{b_{\lambda}}(w,z)dw.
\nonumber
\end{equation}
Thus by Lemma \ref{lemma4}
\begin{eqnarray}
|\partial_{x_i}K^{b_{\lambda}}_{B}(x, z)|
&\le&\int_{B}|\partial_{x_i}G_{B}(x , y)|J^{b_{\lambda}}(y,z)dy \nonumber \\
&&+\int_{B}\left(\int_{B}|\partial_{x_i}G_{B}(x, y) \mathcal{S}^{b_{\lambda}}_{z}G^{b_{\lambda}}_{B}(y,w)|dy\right)
J^{b_{\lambda}}(w,z)dw.
\nonumber
\end{eqnarray}
On the other hand, by \eqref{e:3.3} and \eqref{scalingforgb}, we
have
$$
\frac{1}{2}G_{B}(x,y)\le G^{b_{\lambda}}_{B}(x,y)\le\frac{3}{2}G_{B}(x,y)
\quad\hbox{for }  x, y \in B.
$$
Thus
\begin{equation}
K^{b_{\lambda}}_{B}(x,z)=\int_{B}G^{b_{\lambda}}_{B}(x,y)J^{b_{\lambda}}(y,z)dy\ge
\frac{1}{2}\int_{B}G_{B}(x,y)J^{b_{\lambda}}(y,z)dy
\quad\hbox{for } z\in \bar{B}^{c}.
\nonumber
\end{equation}
Now \eqref{partialK} is  implied by \eqref{e:3.3} and Lemmas \ref{lemma6}-\ref{lemma8}. This completes the proof of
the theorem. \qed

\begin{lemma}\label{lemma9}
Suppose  Assumption 1 holds and $f$ is regular harmonic with respect to
$X^{b}$ in $B(x,r)$ for some $x\in\mathbb{R}^{d}$ and $r\in (0,
r_{1}]$. Then $\partial_{x_i}f(x)$ exists for every $1\leq i\leq d$ and
\begin{equation}\label{partialf}
\partial_{x_i} f(x)=\int_{\overline{B(x, r)}^{c}}f(z)\partial_{x_i}
K^{b}_{B(x, r)}(x, z)dz.
\end{equation}
\end{lemma}

\proof Recall that $e_i$ is the unit vector along the positive $x_i$-axis.
 Choose $\varepsilon>0$ sufficiently small so that
$x+\varepsilon e_{i}\in B(x, r/4)$. By the regular harmonicity of
$f$, we have
$$
\frac{f(x+\varepsilon e_{i})-f(x)}{\varepsilon}=\int_{\overline{B(x,r)}^{c}}f(z)
\left[\frac{K^{b}_{B(x,r)}(x+\varepsilon e_{i},z)-K^{b}_{B(x,r)}(x,z)}{\varepsilon}\right] dz.
$$
Therefore \eqref{partialf} follows from \eqref{partialk} and the
dominated convergence theorem. \qed

\bigskip

\noindent\textbf{Proof of Theorem \ref{theorem1}:} Let $x\in D$ and
$0<r<(\delta_{D}(x)\wedge r_{1})/2$. Note that under our assumption,
$f$ is regular harmonic in $B(x,r)$ with respect to $X^{b}$. By
\eqref{partialf} and \eqref{partialk}, we have
\begin{eqnarray}
|\partial_{x_i} f(x)|&\le&\int_{\overline{B(x,r)}^{c}}f(z)|\partial_{x_i}
K^{b}_{B(x,r)}(x,z)|dz\nonumber\\
&\le&\frac{C_{14}}{r}\int_{\overline{B(x, r)}^{c}}f(z) K^{b}_{B(x, r)}(x,z)dz\nonumber\\
&=&\frac{C_{14}}{r}f(x)
 \rightarrow \frac{2C_{14}}{r_{1}\wedge \delta_{D}(x)}f(x)
\quad \hbox{as }
r\uparrow (r_{1}\wedge \delta_{D}(x))/2.\nonumber
\end{eqnarray}
\qed

\section{Gradient lower bound estimate}\label{S:5}

For   $x=(x_{1},\cdots,x_{d})\in \mathbb{R}^{d}$, we write
$x=(\tilde{x},x_{d})$, where $\tilde{x}=(x_{1},\cdots,x_{d-1})$.
In this section, we fix a Lipschitz function
$\Gamma:\mathbb{R}^{d-1}\to\mathbb{R}$ with Lipschitz constant
$\lambda_{0}$ so that $|\Gamma(\tilde{x})-\Gamma(\tilde{y})|\le
\lambda_{0}|\tilde{x}-\tilde{y}|$ for all
$\tilde{x},\tilde{y}\in\mathbb{R}^{d-1}$. Put
$\rho(x):=x_{d}-\Gamma(\tilde{x})$. Unless stated otherwise, $D$
denotes the special Lipschitz open set defined by
$D=\{x\in\mathbb{R}^{d}:\rho(x)>0\}$.  When $x\in D$,  $\rho(x)$ serves
as the vertical distance from $x\in D$ to $\partial D$, and it satisfies
\begin{equation}
\rho(x)/\sqrt{1+\lambda_{0}^{2}}\le
\delta_{D}(x)\le\rho(x)
\quad\hbox{for }  x\in D.\label{**}
\end{equation}
We define the ``box"
${D^+}(x,h,r):=\{y\in\mathbb{R}^{d}:0<\rho(y)<h,\
|\tilde{x}-\tilde{y}|<r\}$, and the ``inverted box" ${D^-}
(x,h,r):=\{y\in\mathbb{R}^{d}:-h<\rho(y)\le 0,\
|\tilde{x}-\tilde{y}|<r\}$, where $x\in \mathbb{R}^{d}$ and $h,\
r>0$.

\begin{lemma}\label{lemma16}
Let $r_1\in (0, 1]$ be the constant in Lemma \ref{lemma2}.
Suppose Assumption 1 holds, $z_0\in\partial D$ and $r\in (0,r_{1}/2]$. Let $A_{z_0}\in D$ be
such that $\rho(A_{z_0})=|A_{z_0}-z_0|=r/2$. Then there exist positive
constants $C_{15}=C_{15}(d,\alpha,\beta,\lambda_{0},{M_1},{M_2})$ and
${\gamma_2}={\gamma_2}(d,\alpha,\beta,\lambda_{0},{M_1},{M_2})$ such that
for every non-negative function $u$ which is harmonic in $D\cap B(z_0,2r)$ and
vanishes in $D^{c}\cap B(z_0,2r)$, we have
$$\frac{u(x)}{u(A_{z_0})}\ge C_{15}\left(\frac{\rho(x)}{\rho(A_{z_0})}\right)^{\alpha-{\gamma_2}}\quad\hbox{for }  x\in D\cap B(z_0,r).$$
\end{lemma}
\proof Note that by Lemma \ref{lemma2} and Assumption 1, we have for
every $x\in\mathbb{R}^{d}$ and $y\in \overline{B(x,r)}^{c}$,
\begin{equation}\label{lem16.1}
K^{b}_{B(x,r)}(x,y)  \stackrel{{M_2}} \asymp \int_{B(x,r)}G_{B(x,r)}(x,z)J^{\eps_0}(z,y)dz
\asymp K_{B(x,r)}^{\eps_0}(x,y).
\end{equation}
Lemma \ref{lemma16} follows from \eqref{lem16.1},
the uniform Harnack inequality (Theorem \ref{HI1})
and a standard argument of induction in the same way as for the case
of symmetric $\alpha$-stable process in
\cite[Lemma 5]{Bogdan1} (see also \cite[Lemma 4.2]{BKN}).
We omit the details here.\qed

In the remaining of this section, we assume that
Assumption 1 and Assumption 2  with $i=d$ hold.
In this case, the jumping kernel $J^{b}(x,y)$ of the
process $X^{b}$ satisfies that for every $x\in\mathbb{R}^{d}$,
\begin{equation}
J^{b}(x,y)=\frac{\mathcal{A}(d,-\alpha)}{|y-x|^{d+\alpha}}+\mathcal{A}(d,-\beta)\frac{\varphi(\tilde{x})\psi(|y-x|)}{|y-x|^{d+\beta}}=:j^{b}(\tilde{x},|y-x|)\quad\mbox{a.e.
}y\in\mathbb{R}^{d}.\label{jb}
\end{equation}
We note that by condition \eqref{e:1.10} of Assumption 2,
$j^{b}(\tilde{x},|z|)$ is non-increasing in $|z|$ for every
$\tilde{x}\in\mathbb{R}^{d-1}$.
Fix $z_0\in \partial D$, $r\in (0,r_{1}]$. We
define
${D^+}:={D^+}(z_0,4r\sqrt{1+\lambda_{0}^{2}},2r)\supset
D\cap B(z_0,2r)$, and
\begin{equation}
g_{b,r}(x):=\mathbb{P}_{x}\left(X^{b}_{\tau_{{D^+}
}}\not\in {D^-}(z_0,\infty,2r)\right)
\quad \hbox{for }  x\in \mathbb{R}^{d}.
\nonumber
\end{equation}
Clearly, $g_{b,r}$ is regular harmonic in ${D^+} $ with
respect to $X^{b}$, $g_{b,r}(x)=0$ in
${D^-}(z_0,\infty,2r)$, and $g_{b,r}(x)=1$ in
$({D^+} \cup{D^-}(z_0,\infty,2r))^{c}$.

\begin{lemma}
The function $g_{b,r}(x)$ is non-decreasing in $x_{d}$.
\end{lemma}
\proof Note that
$g_{b,r}(x)=1-\mathbb{P}_{x}\left(X^{b}_{\tau_{{D^+}
}}\in {D^-}(z_0,\infty,2r)\right)$ for every
$x\in\mathbb{R}^{d}$. Take $x,y\in {D^+}$ such that
$\tilde{x}=\tilde{y}$ and
$y_{d}< x_{d}$.
Consider the
process $(X^{b}_{t},\mathbb{P}_{x})$
starting from $x$
(i.e. $\mathbb{P}_{x}(X^{b}_{0}=x)=1$). For every $t\ge 0$, define
$$Y^{b}_{t}:=X^{b}_{t}- (x_d-y_d) e_d.$$
Then $(Y^{b}_{t},\mathbb{P}_{x})$ is a Markov process
starting from $y$.
Let $\mathcal{S}(\mathbb{R}^{d})$ denote the totality of tempered functions on $\mathbb{R}^{d}$. For every $f\in\mathcal{S}(\mathbb{R}^{d})$, if we define $f_{d}(z):=f(z-(x_d-y_d) e_d)$ for $z\in\mathbb{R}^{d}$, then
$$\mathcal{L}^{b}f(x-(x_{d}-y_{d})e_{d})=\mathcal{L}^{b}f_{d}(x)=\lim_{\varepsilon\to 0}\int_{|y-x|>\varepsilon}(f_{d}(y)-f_{d}(x))j^{b}(\tilde{x},|y-x|)dy\mbox{ for }x\in\mathbb{R}^{d}.$$
Thus
$$
f(Y^{b}_{t})-f(Y^{b}_{0})-\int_{0}^{t}\mathcal{L}^{b}f(Y^{b}_{s})ds
=f_{d}(X^{b}_{t})-f_{d}(X^{b}_{0})-\int_{0}^{t}\mathcal{L}^{b}f_{d}(X^{b}_{s})ds
$$
is a $\mathbb{P}_{x}$-martingale.
We know from
\cite[Theorem 1.3]{CW}
that the solution of the martingale problem $(\mathcal{L}^{b},\mathcal{S}(\mathbb{R}^{d}))$ with initial value $y$ is unique.
Hence $(Y^{b}_{t},\mathbb{P}_{x})$ has the
same distribution as $(X^{b}_{t},\mathbb{P}_{y})$. Consider the
trajectory $\omega$ of $X^{b}_{t}$
starting from $x$.
If $\omega$
exits ${D^+}$ by going into
${D^-}(0,+\infty,2r)$, then so does $\omega-(x_d-y_d) e_d$ which
is the trajectory of $Y^{b}_{t}$
starting from $y$.  Hence
$$
\mathbb{P}_{x}\left(X^{b}_{\tau_{{D^+} }}\in
{D^-}(z_0,\infty,2r)\right)\le
\mathbb{P}_{x}\left(Y^{b}_{\tau_{{D^+} }}\in
{D^-}(z_0,\infty,2r)\right)=
\mathbb{P}_{y}\left(X^{b}_{\tau_{{D^+} }}\in
{D^-}(z_0,\infty,2r)\right).
$$
This completes the proof.       \qed

\begin{lemma}\label{lemma17}
Let $r_1\in (0, 1]$ be the constant in Lemma \ref{lemma2}.
There are constants
\hfill \break
$C_{16}=C_{16}(d,\alpha,\beta,\lambda_{0},{M_1},{M_2})>0$
and ${r_2}={r_2}(d, \alpha, \beta, \lambda_{0},{M_1},{M_2}) \in (0, r_1]$
such that for every $z_0\in \partial D$
and $r\in (0, r_2]$,
\begin{equation}
\partial_{x_d}g_{b,r}(x)\ge
C_{16}\frac{g_{b,r}(x)}{\delta_{D}(x)}\quad\hbox{for }  x\in D\cap
B(z_0,r/2).
\nonumber
\end{equation}
\end{lemma}

\proof Without loss of generality we assume $z_0=0$.
Let $r_2 \in (0, r_1]$ to be specified later.
For $r\in (0, r_2]$,  fix $x\in D\cap B(0,r/2)$.
By \eqref{**},
$r_{0}:=\rho(x)/2\sqrt{1+\lambda_{0}^{2}}\le \delta_{D}(x)/2\le r/4\leq r_2/4$.
Set
$${\widehat x}=x+2\rho(x)e_{d} \quad \hbox{and} \quad  {\check{x}}=x-2\rho(x)e_{d}.$$
Observe that
$B(x, r_0),B({\widehat x}, r_0)\subset{D^+}$ and
$B({\check{x}}, r_0)\subset {D^-}(0,+\infty,2r)$. By
\eqref{partialf} and \eqref{2.21}, we have
\begin{eqnarray}
 \partial_{x_d} g_{b,r}(x)
&=&\int_{\overline{B(x, r_0)}^{c}}g_{b,r}(z)\partial_{x_d}
K^{b}_{B(x, r_0)}(x,z)dz\nonumber\\
&\ge&\int_{\overline{B(x, r_0)}^{c}}g_{b,r}(z)\left[\int_{B(x, r_0)}\partial_{x_d}G_{B(x, r_0)}(x,y)J^{b}(y,z)dy\right]dz  \label{3.17} \\
&&-\int_{\overline{B(x, r_0)}^{c}}g_{b,r}(z)\left(\int_{B(x, r_0)\times B(x, r_0)}\left|\partial_{x_d}G_{B(x, r_0)}(x,y)S^{b}_{y}G^{b}_{B(x, r_0)}(y,w)\right|J^{b}(w,z)dy\,dw\right)dz . \nonumber
\end{eqnarray}
Let $\lambda:=1/r_{0}$ and $B_{1}:=\lambda B(x, r_0)$. By
scaling property, Lemma \ref{lemma8} and Lemma \ref{lemma2}, we have
\begin{eqnarray}
&&\int_{B(x, r_0)}\int_{B(x, r_0)}\left|\partial_{x_d}G_{B(x, r_0)}(x,y)S^{b}_{y}G^{b}_{B(x, r_0)}(y,w)\right|J^{b}(w,z)dy\,dw\nonumber\\
&=&\lambda^{d+1}\int_{B_{1}}\int_{B_{1}}\left|\partial_{x_d}G_{B_{1}}(\lambda x,y)S^{b_{\lambda}}_{y}G^{b_{\lambda}}_{B_{1}}(y,w)\right|J^{b_{\lambda}}(w,\lambda z)dy\,dw\nonumber\\
&\le&c_{1}\lambda^{d+1-\alpha+\beta}\int_{B_{1}}G_{B_{1}}(\lambda
x,w)J^{b_{\lambda}}(w,\lambda z)dw\nonumber\\
&=&c_{1}\lambda^{1-\alpha+\beta}\int_{B(x, r_0)}G_{B(x, r_0)}(x,v)J^{b}(v,z)dv\nonumber\\
&\le&2c_{1}\lambda^{1-\alpha+\beta}\int_{B(x, r_0)}G^{b}_{B(x, r_0)}(x,v)J^{b}(v,z)dv\nonumber\\
&=&2c_{1}r_{0}^{-1+\alpha-\beta}K^{b}_{B(x, r_0)}(x,z).
\nonumber
\end{eqnarray}
Here $c_{1}=c_{1}(d,\alpha,\beta,{M_1},{M_2})>0$.
Thus we can continue the estimate in \eqref{3.17} to get
\begin{eqnarray}
&&\partial_{x_d} g_{b,r}(x)\nonumber\\
&\ge&\int_{\overline{B(x, r_0)}^{c}}g_{b,r}(z)\left(\int_{B(x, r_0)}\partial_{x_d}G_{B(x, r_0)}(x,y)J^{b}(y,z)dy\right)dz\nonumber\\
&&-\frac{2c_{1}}{r_{0}}r_{0}^{\alpha-\beta}\int_{\overline{B(x, r_0)}^{c}}g_{b,r}(z)K^{b}_{B(x, r_0)}(x,z)dz\nonumber\\
&=&\int_{\overline{B(x, r_0)}^{c}}g_{b,r}(z)\left(\int_{B(x, r_0)}\partial_{x_d}G_{B(x, r_0)}(x,y)J^{b}(y,z)dy\right)dz-\frac{2c_{1}}{r_{0}}r_{0}^{\alpha-\beta}g_{b,r}(x).\label{partialg}
\end{eqnarray}
Note that by \eqref{ballgreen},
\begin{eqnarray} \label{partialballgreen}
&&\partial_{x_d}G_{B(x, r_0)}(x,y)\nonumber\\
&=&2^{1-\alpha}\pi^{-d/2}\Gamma\left(\frac{d}{2}\right)\Gamma\left(\frac{\alpha}{2}\right)^{-2}
\frac{y_{d}-x_{d}}{|y-x|^{2}}r_{0}^{\alpha}
(r_{0}^{2}-|y-x|^{2})^{\alpha/2}\left(r_{0}^{2}(r_{0}^{2}-|y-x|^{2})+|y-x|^{2}\right)^{-d/2}
\nonumber\\
&&+(d-\alpha)\frac{y_{d}-x_{d}}{|y-x|^{2}}G_{B(x, r_0)}(x,y).
\end{eqnarray}
Obviously $\partial_{x_d}G_{B(x, r_0)}(x,y)$ is anti-symmetric in $y$
with respect to the hyperplane $\mathcal{H}:=\{y\in\mathbb{R}^{d}:\
y_{d}=x_{d}\}$. For every $z\in\overline{B(x, r_0)}^{c}$,
define
$h_{x}(z):=\int_{B(x, r_0)}\partial_{x_d}G_{B(x, r_0)}(x,y)J^{b}(y,z)dy.$
By \eqref{jb} we have for a.e.
$z\in\overline{B(x, r_0)}^{c}$,
\begin{eqnarray}
h_{x}(z)&=&\int_{\{y\in
B(x, r_0),y_{d}>x_{d}\}}\partial_{x_d}G_{B(x, r_0)}(x,y)\left(J^{b}(y,z)-J^{b}(\bar{y},z)\right)dy\nonumber\\
&=&\int_{\{y\in
B(x, r_0),y_{d}>x_{d}\}}\partial_{x_d}G_{B(x, r_0)}(x,y)\left(j^{b}(\tilde{y},|z-y|)-j^{b}(\tilde{y},|z-\bar{y}|)\right)dy,
\label{3.21}
\end{eqnarray}
where $\bar{y}=(\tilde{y},2x_{d}-y_{d})$. We observe that the right
hand side of \eqref{3.21} is antisymmetric in $z$ with respect to
the hyperplane $\mathcal{H}$. Recall that $j^{b}(\tilde{y},r)$ is
non-increasing in $r$. Following from this, the monotonicity of
$g_{b,r}(x)$ in $x_{d}$ and the Harnack inequality, we have
\begin{eqnarray}
 \int_{\overline{B(x, r_0)}^{c}}g_{b,r}(z)h_{x}(z)dz
&\ge& \int_{B({\widehat x}, r_0)\cup B({\check{x}}, r_0)}g_{b,r}(z)h_{x}(z)dz
 =
\int_{B({\widehat x}, r_0)}g_{b,r}(z)h_{x}(z)dz    \nonumber\\
&\ge& c_{2}g_{b,r}(x) \int_{B({\widehat x}, r_0)}h_{x}(z)dz
\label{18}
\end{eqnarray}
for some $c_{2}=c_{2}(d,\alpha,\beta,\lambda_{0},{M_1},{M_2})>0$.
Note that by \eqref{partialballgreen} and \eqref{3.21}, for a.e.
$z\in B({\widehat x}, r_0)$,
\begin{eqnarray}
h_{x}(z)&\ge&(d-\alpha)\int_{\{y\in
B(x, r_0):y_{d}>x_{d}\}}\frac{y_{d}-x_{d}}{|y-x|^{2}}G_{B(x, r_0)}(x,y)\left(j^{b}(\tilde{y},|z-y|)-j^{b}(\tilde{y},|z-\hat{y}|)\right)dy\nonumber\\
&=&(d-\alpha)\int_{B(x, r_0)}\frac{y_{d}-x_{d}}{|y-x|^{2}}G_{B(x, r_0)}(x,y)j^{b}(\tilde{y},|z-y|)dy>0.\nonumber
\end{eqnarray}
Therefore  by the scaling property of the Green function $G_{B(x, r_0)}$,
\begin{eqnarray}
&&\int_{B({\widehat x}, r_0)}h_{x}(z)dz\nonumber\\
&\ge&(d-\alpha)\int_{B({\widehat x}, r_0)}\int_{B(x, r_0)}\frac{y_{d}-x_{d}}{|y-x|^{2}}G_{B(x, r_0)}(x,y)J^{b}(y,z)dy\,dz\nonumber\\
&\ge& \frac{d-\alpha}{M_2} \int_{|v-\frac{{\widehat x}}{r_{0}}|<1}\int_{|w-\frac{x}{r_{0}}|<1}\frac{1}{r_{0}}\frac{w_{d}-x_{d}/r_{0}}{|w-x/r_{0}|^{2}}G_{B(0,r_{0})}(0,r_{0}(w-x/r_{0}))J(r_{0}|v-w)|)r_{0}^{2d}dw\,dv\nonumber\\
&=&\frac{d-\alpha}{r_{0}M_2}\int_{|v-\frac{{\widehat x}}{r_{0}}|<1}\int_{|w-\frac{x}{r_{0}}|<1}\frac{w_{d}-x_{d}/r_{0}}{|w-x/r_{0}|^{2}}G_{B(0,1)}(0,w-x/r_{0})J(|v-w|)dw\,dv\nonumber\\
&=&\frac{d-\alpha}{r_{0} M_2}\int_{|v|<1}\int_{|w|<1}\frac{w_{d}}{|w|^{2}}G_{B(0,1)}(0,w)J(|v-w+\frac{{\widehat x}-x}{r_{0}}|)dw\,dv\nonumber\\
&=&\frac{d-\alpha}{r_{0}M_2}\int_{|v|<1}\int_{|w|<1}\frac{w_{d}}{|w|^{2}}G_{B(0,1)}(0,w)J(|v-w
+4\sqrt{1+\lambda_{0}^{2}}e_d|)dw\,dv=:\frac{c_3}{r_{0}}, \label{17}
\end{eqnarray}
with $c_{3}=c_{3}(d,\alpha,{M_2})>0$.
It follows from \eqref{partialg}, \eqref{partialballgreen}, \eqref{18} and \eqref{17} that
\begin{equation}
\partial_{x_d}g_{b,r}(x)\ge
\frac{1}{r_{0}}\left( c_{2}c_{3}-2c_{1}r_{0}^{\alpha-\beta}\right) g_{b,r}(x)\ge
\frac{2}{\delta_{D}(x)}\left( c_{2}c_{3}-2c_{1}\left(r_2/4\right)^{\alpha-\beta}\right) g_{b,r}(x).\label{3.24}
\end{equation}
The lemma now follows from \eqref{3.24} by setting $r_2$ so
small that $2c_{1}\left({r_2}/{4}\right)^{\alpha-\beta}\le
c_{2}c_{3}/2$.\qed

\bigskip

\begin{lemma}\label{lemma18}
Let
$r_2\in (0, r_1]\subset (0, 1]$ be the constant in Lemma \ref{lemma17}.
There is a positive constant $C_{17}=C_{17}(d,\alpha,\beta,\lambda_{0},{M_1},{M_2})$
such that for every $r\in (0,r_{2}]$,
there is a constant $r_{3}=r_{3}(d,\alpha,\beta,\lambda_{0},M_{1},M_2,r)\in (0,r/2)$ so that for every $z_0\in\partial D$  and every non-negative function $f$
that is regular harmonic in $D\cap B(z_0,2r)$ with
respect to $X^{b}$ and vanishes in $D^{c}\cap B(z_0,2r)$,
$$\left|\partial_{x_d}f(x)\right|\ge C_{17}\frac{f(x)}{\delta_{D}(x)}\quad \hbox{for }
x\in D\cap B(z_{0},r_{3}).
$$
\end{lemma}

\proof   Without loss of generality we assume $z_0=0$.
For $r\in (0,r_{2}]$, fix an arbitrary $x\in D\cap B(0,r/(2\sqrt{1+\lambda_{0}^{2}}))$.
Let $z_{x}\in\partial D$ be such that $|x-z_{x}|=\rho(x)$. Define $c:=\lim_{D\ni
y\to z_{x}}f(y)/g_{b,r}(y)$ and $u(y):=c\,g_{b,r}(y)$. Obviously
$B(z_{x},3r/2)\subset B(0,2r)$, and thus $f$, $u$ are harmonic in
$D\cap B(z_{x},3r/2)$ and vanish in $D^{c}\cap B(z_{x},3r/2)$.
Since $\lim_{D\ni y\to z_{x}}u(y)/f(y)=\lim_{D\ni y\to
z_{x}}f(y)/u(y)=1$, by Lemma \ref{lemma14}, for any $y\in D\cap
B(z_{x},3r/4)$,
\begin{equation}
\left|\frac{u(y)}{f(y)}-1\right|\vee\left|\frac{f(y)}{u(y)}-1\right|\le
c_{1}\left(\frac{|y-z_{x}|}{r}\right)^{{\gamma_1}}\le
c_{2}\label{eq10}
\end{equation}
for some positive constants $c_{i}=c_{i}(d,\alpha,\beta,\lambda_{0},{M_1},{M_2})$, $i=1,2$.
Consequently,
\begin{equation}
(1+c_{2})^{-1}f(y)\le u(y)\le (1+c_{2})f(y)\quad\hbox{for }  y\in D\cap
B(z_{x},3r/4).\label{eq11}
\end{equation}
Note that $x\in D\cap B(z_{x},3r/4)$ since $\rho(x)\le \delta_{D}(x)\sqrt{1+\lambda_{0}^{2}}\le r/2$.
By Lemma \ref{lemma17} and \eqref{eq11}, we have
\begin{eqnarray}
\partial_{x_d}f(x)&\ge& \partial_{x_d}u(x)-\left|\partial_{x_d}(f-u)(x)\right|
\geq
c_{3}\frac{u(x)}{\delta_{D}(x)}-\left|\partial_{x_d}(f-u)(x)\right| \nonumber\\
&\ge& c_{4}\frac{f(x)}{\delta_{D}(x)}-\left|\partial_{x_d}(f-u)(x)\right|.     \label{eq6}
\end{eqnarray}
We assume $\rho(x)<3r/128$.
Set $v(y):=f(y)-u(y)$ and $\xi:=2\rho(x)$.
Let $\eta\in (16\rho(x), 3r/8)$ to be specified later.
In the rest of this proof, we set
${D_1}:={D^+}(z_{x},\xi,\xi)$ and
${D_2}:={D^+}(z_{x},\eta,\eta)$. Then
$B(x,\delta_{D}(x))\subset{D_1}\subset {D_2}\subset D\cap
B(z_{x},3r/4)$ and $\delta_{D}(x)=\delta_{{D_1}}(x)$.
Define
$V(y):=\E_{y}\left[|v|(X^{b}_{\tau_{{D_1}}})\right]$.
Clearly $V$ is regular harmonic in ${D_1}$ with respect to
$X^{b}$ and $|v(y)|\le V(y)$ for all $y\in\mathbb{R}^{d}$. By
Theorem \ref{theorem1}, we have
\begin{eqnarray}
|\partial_{x_d}v(x)|&\le&|\partial_{x_d}V(x)|+|\partial_{x_d}(V-v)(x)|
\le c_{5}\frac{V(x)}{\delta_{{D_1}}(x)}=c_{5}\frac{V(x)}{\delta_{D}(x)}.\label{eq7}
\end{eqnarray}
We aim to estimate $V(x)$. Note that
\begin{eqnarray}
V(x)&\le&\mathbb{E}_{x}\left[|v|(X^{b}_{\tau_{{D_1}}}):X^{b}_{\tau_{{D_1}}}\in
{D_2}\right]+\mathbb{E}_{x}\left[f(X^{b}_{\tau_{{D_1}}}):X^{b}_{\tau_{{D_1}}}\in
D^c_2\right]+\mathbb{E}_{x}\left[u(X^{b}_{\tau_{{D_1}}}):X^{b}_{\tau_{{D_1}}}\in
D^c_2\right]\nonumber\\
&=:&I(x)+II(x)+III(x).\nonumber
\end{eqnarray}
By \eqref{eq10}, for any $y\in {D_2}\subset D\cap
B(z_{x},2\eta)\subset D\cap B(z_{x},3r/4)$, we have
$$|v(y)|=f(y)\left|\frac{u(y)}{f(y)}-1\right|\le c_{1}f(y)\left(\frac{|y-z_{x}|}{r}\right)^{{\gamma_1}}\le c_{6}\left(\frac{\eta}{r}\right)^{{\gamma_1}}f(y).$$ Thus
\begin{equation}
I(x)\le
c_{6}\left(\frac{\eta}{r}\right)^{{\gamma_1}}\mathbb{E}_{x}\left[f(X^{b}_{\tau_{{D_1}}})\right]=c_{6}\left(\frac{\eta}{r}\right)^{{\gamma_1}}f(x).\label{I}
\end{equation}
Let $A_{x}\in D$ be such that $\rho(A_{x})=|A_{x}-z_{x}|=\eta/16$.
Define ${D_3}:=B(A_{x},\eta/16\sqrt{1+\lambda_{0}^{2}})$. We
observe that ${D_3}\subset D\cap B(z_{x},\eta/2)\subset {D_2}$ and
${D_1}\subset D\cap B(z_{x},\eta/4)$. For any $y\in
D^{c}_{2}\cap \mbox{supp}f$ and $z\in {D_1}$, we have
$|y-A_{x}|\ge 7\eta/16$, $|A_{x}-z|\le 5\eta/16$, and
\begin{equation}
|y-z|\ge |y-A_{x}|-|A_{x}-z|\ge \frac{2}{7}|y-A_{x}|.\label{3.39}
\end{equation}
If we let $\lambda:=1/\mathrm{diam}({D_1})$, then
$\|b_{\lambda}\|_{\infty}=\mathrm{diam}({D_1})^{\alpha-\beta}\|b\|_{\infty}\le
(8\rho(x))^{\alpha-\beta}{M_1}\le {M_1}$.
Thus by \eqref{scalingforgb} and
Lemma \ref{lemma10}, we have
\begin{equation}\label{3.40}
G^{b}_{{D_1}}(x,z)=\lambda^{d-\alpha}G^{b_{\lambda}}_{\lambda
D}(\lambda x,\lambda z)\le c_{7}|x-z|^{\alpha -d},
\quad  z\in{D_1}
\end{equation}
for some constant $c_{7}=c_{7}(d,\alpha,\beta,{M_1})>0$.
So by \eqref{3.39} and \eqref{3.40}, for any $y\in
D^{c}_{2}\cap \mbox{supp}f$,
\begin{eqnarray}
K^{b}_{{D_1}}(x,y)&=&\int_{{D_1}}G^{b}_{{D_1}}(x,z)J^{b}(z,y)dz
\leq
c_{7}{M_2}\int_{{D_1}}|x-z|^{\alpha -d}J^{\eps_0}(|y-z|)dz\nonumber\\
&\lesssim&c_{7}{M_2}
J^{\eps_0}(|y-A_{x}|)\int_{B(z_{x},2\xi)}|x-z|^{\alpha -d}dz
\lesssim
c_{7}{M_2}\xi^{\alpha}J^{\eps_0}(|y-A_{x}|).\label{eq8}
\end{eqnarray}
On the other hand for any
$y\in D^{c}_{2}\cap \mbox{supp}f$
and
$z\in{D_3}$, we have $|y-z|\le |y-A_{x}|+|A_{x}-z|\le
12|y-A_{x}|/7$. Thus by Lemma \ref{lemma2} and \eqref{condi2}
\begin{eqnarray}
K^{b}_{{D_3}}(A_{x},y)&=&\int_{{D_3}}G^{b}_{{D_3}}(A_{x},z)J^{b}(z,y)dz\nonumber\\
&\gtrsim&M^{-1}_{2}\left(\int_{{D_3}}G_{{D_3}}(A_{x},z)dz\right)J^{\eps_0}(|y-A_{x}|)\nonumber\\
&\asymp&M^{-1}_{2}\eta^{\alpha}J^{\eps_0}(|y-A_{x}|).\label{eq9}
\end{eqnarray}
Combining \eqref{eq8} and \eqref{eq9}, we have
\begin{equation}
K^{b}_{{D_1}}(x,y)\lesssim \frac{\xi^{\alpha}}{\eta^{\alpha}}\,K^{b}_{{D_3}}(A_{x},y)
\quad \hbox{for }   y\in D^{c}_{2}\cap \mbox{supp}f. \label{3.32}
\end{equation}
Consequently, by \eqref{3.32}, Lemma \ref{lemma16} and \eqref{**},
we have
\begin{eqnarray}
II(x)&=&\int_{D^{c}_{2}\cap
\mbox{supp}f}f(y)K^{b}_{{D_1}}(x,y)dy
 \lesssim \frac{\xi^{\alpha}}{\eta^{\alpha}}\int_{D^{c}_{2}\cap
\mbox{supp}f}f(y)K^{b}_{{D_3}}(A_{x},y)dy\nonumber\\
&\le&\frac{\xi^{\alpha}}{\eta^{\alpha}}\int_{D^{c}_{3}}f(y)K^{b}_{{D_3}}(A_{x},y)dy
 \, = \, \frac{\xi^{\alpha}}{\eta^{\alpha}}f(A_{x}) \nonumber \\
&\lesssim&\frac{\xi^{\alpha}}{\eta^{\alpha}}\,\frac{\rho(A_{x})^{\alpha-{\gamma_2}}}
{\rho(x)^{\alpha-{\gamma_2}}}\,f(x)
\, \asymp \, \frac{\rho(x)^{{\gamma_2}}}{\eta^{{\gamma_2}}}\,f(x).\label{II}
\end{eqnarray}
Similarly we can prove that
\begin{equation}
III(x)\lesssim
\frac{\rho(x)^{{\gamma_2}}}{\eta^{{\gamma_2}}}\,u(x)\lesssim
\frac{\rho(x)^{{\gamma_2}}}{\eta^{{\gamma_2}}}\,f(x).\label{III}
\end{equation}
Combining \eqref{I}, \eqref{II} and \eqref{III}, we have
\begin{equation} \label{3.46}
V(x)\le\left( c_{6}\left(\frac{\eta}{r}\right)^{{\gamma_1}}+c_{8}
\left(\frac{\rho(x)}{\eta}\right)^{{\gamma_2}}\right) f(x).
\end{equation}
Thus by \eqref{eq6}, \eqref{eq7} and \eqref{3.46} we have
\begin{equation} \label{eq12}
\partial_{x_d}f(x)\ge
\left(c_{4}-c_{5}\left(c_{6}\frac{\eta^{{\gamma_1}}}
{r^{{\gamma_1}}}+c_{8}\frac{\rho(x)^{{\gamma_2}}}{\eta^{{\gamma_2}}}\right)\right)
\frac{f(x)}{\delta_{D}(x)}.
\end{equation}
Let $\eta=16\rho(x)^{{\gamma_2}/({\gamma_1}+{\gamma_2})}$.
The lemma now follows
from \eqref{eq12} and \eqref{**} provided we choose
$r_{3}$ small enough such that
$c_{5}\left(c_{6}16^{{\gamma_1}}r^{-{\gamma_1}}+c_{8}16^{-{\gamma_2}}\right)
(r_3 (1+\lambda_0))^{{{\gamma_1}{\gamma_2}}/{(\gamma_1+\gamma_2)}}\le c_{4}/2$. \qed
\bigskip

\noindent \textbf{Proof of Theorem \ref{theorem2}:} The upper bound in
\eqref{them2.main} was
established
 more generally in Theorem
\ref{theorem1}. The lower bound follows from Lemma \ref{lemma18} and
the inequality $|\nabla f|\ge |\partial_{x_d}f|$.\qed

\bigskip

\begin{remark}\rm
Taking $b(x, z)=\eps\in (0, M_2]$ in Assumption 1, we get from
Theorem \ref{theorem1} and Theorem \ref{T:1.7}
the uniform gradient estimate for mixed stable processes, which
in particular recovers a main result of \cite{Y} on gradient estimates.
\end{remark}

\bigskip

\section{Examples}

In this section,
we give some concrete examples where
Assumptions 1 and 3 hold.

\begin{example} \rm
If $b(x,z)=1_{\{|z|\le c_{1}\}}$ for some $c_{1}>0$, the jumping kernel of the corresponding Feller process $X^{b}$ is
$$J^{b}(x,y)=\frac{\mathcal{A}(d,-\alpha)}{|x-y|^{d+\alpha}}+\frac{\mathcal{A}(d,-\beta)}{|x-y|^{d+\beta}}\,1_{\{|x-y|\le c_{1}\}}.$$
In this case $X^{b}$ is the independent sum of a symmetric $\alpha$-stable process and a truncated symmetric $\beta$-stable process, and
Assumptions 1 and 3 hold with $\eps_0=0$ and $\psi(r)=1_{\{r\le c_{1}\}}$,
respectively.

More generally, suppose $b(x,z)=b_1 (x, z)1_{\{|z|\le c_{1}\}}$ for some
 $c_1>0$ and a bounded function $b_1(x, z)$ on $\mathbb{R}^{d}\times \mathbb{R}^{d}$ that is symmetric in $z$ and is bounded between two positive constants.
Then  Assumption 1 holds with $\eps_0=0$. \qed
\end{example}

\begin{example} \rm
If $b(x,z)=1+\frac{\mathcal{A}(d,-\gamma)}{\mathcal{A}(d,-\beta)}|z|^{\beta-\gamma}1_{\{|z|\le c_{2}\}}$ for some $c_{2}>0$ and $0<\gamma<\beta$, the jumping kernel of the corresponding Feller process $X^{b}$ is
$$J^{b}(x,y)=\frac{\mathcal{A}(d,-\alpha)}{|x-y|^{d+\alpha}}+\frac{\mathcal{A}(d,-\beta)}{|x-y|^{d+\beta}}
+\frac{\mathcal{A}(d,-\gamma)}{|x-y|^{d+\gamma}}\,1_{\{|x-y|\le c_{2}\}}.$$
In this case $X^{b}$ is the independent sum of a mixed-stable process and a truncated symmetric $\gamma$-stable process, and
Assumptions 1 and 3 hold with $\eps_0=1$ and
$\psi(r)=1+\frac{\mathcal{A}(d,-\gamma)}{\mathcal{A}(d,-\beta)}\,r^{\beta-\gamma}1_{\{r\le c_{2}\}}$, respectively.

More generally, suppose $b(x,z)$ is a bounded function on $\mathbb{R}^{d}\times \mathbb{R}^{d}$
that is symmetric in $z$ and is bounded between two positive constants.
Then  Assumption 1 holds with $\eps_0=1$.  \qed
\end{example}

\begin{example}\rm
We consider the following stochastic differential equation on $\mathbb{R}^{d}$:
\begin{equation}
dX_{t}=dY_{t}+C(X_{t-})dZ_{t},\label{sde}
\end{equation}
where $Y$ is a symmetric $\alpha$-stable process, $Z$ is an independent $\beta$-stable process with $0<\beta<\alpha$,
and $C$ is a bounded Lipschitz function on $\mathbb{R}^{d}$. Using Picard's iteration method, one can show that for every $x\in\mathbb{R}^{d}$, SDE \eqref{sde} has a unique strong solution with $X_{0}=x$. The collection of the solutions $(X_{t},\mathbb{P}_{x},x\in\mathbb{R}^{d})$ forms a strong Markov process $X$ on $\mathbb{R}^{d}$. Using Ito's formula, one concludes that the infinitesimal generator of $X$ is $\mathcal{L}^{b}$ with $b(x,z)=|C(x)|^{\beta}$. If there exists $c_{3}>0$ such that $|C(x)|\ge c_{3}$ for $x\in\mathbb{R}^{d}$, then our Assumption 1 holds
with $\eps_0=1$. \qed
\end{example}

\bigskip

\bigskip
\noindent\textbf{Acknowledgement} The research of Zhen-Qing Chen is partially supported
by NSF grant DMS-1206276 and NNSFC 11128101. The research of Yan-Xia Ren is supported by NNSFC (Grant No.  11271030 and 11128101). The research of Ting Yang is partially supported by NNSF of China (Grant No. 11501029) and Beijing Institute of Technology Research Fund Program for Young Scholars.

\end{document}